\documentclass[a4paper,10pt]{amsart}

\usepackage[latin1]{inputenc}
\usepackage[american]{babel}
\usepackage{amsmath}
\usepackage{amsfonts}
\usepackage{amsthm}
\usepackage{amssymb}
\numberwithin{equation}{section}

\newcommand{\de}{{\delta}}
\newcommand{\Hn}{\mathbb H^n}
\renewcommand{\H}{\mathbb H}
\newcommand{\Rn}{\mathbb R^n}
\newcommand{\Rm}{\mathbb R^m}
\newcommand{\norm}[1]{\|{#1}\|_{\infty}}
\newcommand{\dede}[1]{{\frac{\partial}{\partial #1}}}
\newcommand{\eps}{{\epsilon}}
\newcommand{\spt}{{\mathrm{spt}}}

\newcommand{\rn}[1]{{\mathbb R}^{#1}}

\newcommand{\norma}{\vert\!\vert}
\renewcommand{\rho}{\varrho}

\renewcommand{\t}{{\tau}}
\newcommand{\ep}{{\epsilon}}
\renewcommand{\d}{{\delta}}
\newcommand{\e}{{\eta}}
\renewcommand{\a}{{\alpha}}
\renewcommand{\l}{{\lambda}}
\renewcommand{\i}{{\iota}}
\newcommand{\z}{{\zeta}}
\newcommand{\G}{{\mathbb G}}

\newcommand{\f}{{\phi}}

\newcommand{\F}{{\Phi}}
\newcommand{\ci}{{\mathbf C}}
\newcommand{\heis}{{\mathbb{H}}}
\newcommand{\p}{{\psi}}
\newcommand{\R}{\mathbb R}

\newcommand{\Xt}{\tilde X}
\newcommand{\Yt}{\tilde Y}
\newcommand{\Tt}{\tilde T}
\newcommand{\bur}{\Wf}
\newcommand{\Wfs}{W^{\f_s}}
\renewcommand{\b}{\beta}
\newcommand{\res}{\mathop{\hbox{\vrule height 7pt width .5pt depth 0pt
\vrule height .5pt width 6pt depth 0pt}}\nolimits}

\newcommand{\average}{{\mathchoice {\kern1ex\vcenter{\hrule height.4pt
width 6pt
depth0pt} \kern-9.7pt} {\kern1ex\vcenter{\hrule height.4pt width 4.3pt
depth0pt}
\kern-7pt} {} {} }}

\newcommand{\N}{\mathbb N}
\newcommand{\Wf}{W^\f}
\newcommand{\vf}{\varphi}
\renewcommand{\div}{\mathrm{div}}

\begin{document}

\newtheorem{theo}{Theorem}[section]
\newtheorem{defi}[theo]{Definition}
\newtheorem{lemma}[theo]{Lemma}
\newtheorem{prop}[theo]{Proposition}
\newtheorem{obs}[theo]{Remark}
\newtheorem{cor}[theo]{Corollary}
\newtheorem{ex}[theo]{Example}

\makeatletter
\def\@eqnnum{\hbox to .01\p@{}\rlap{\reset@font\rm
        \hskip -13.4cm(\theequation)}}
\makeatother


\title[The Bernstein problem for intrinsic graphs in $\Hn$ and calibrations]{The Bernstein problem for intrinsic graphs in Heisenberg groups and calibrations }
\author{Vittorio Barone Adesi}
\address{Vittorio Barone Adesi: Dipartimento di Matematica\\Universit\`a di Trento\\ Via Sommarive 14\\ 38050, Povo (Trento) - Italy\\} 
\email{vbarone@science.unitn.it}
\thanks{V.B.A. is supported by MIUR, Italy, GNAMPA of INDAM and University of Trento, Italy.}
\author{Francesco Serra Cassano}
\address{Francesco Serra Cassano: Dipartimento di Matematica\\Universit\`a di Trento\\ Via Sommarive 14\\ 38050, Povo (Trento) - Italy\\} 
\email{cassano@science.unitn.it}
\thanks{F.S.C. is supported by MIUR, Italy, GNAMPA of INDAM and University of Trento, Italy.}
\author{Davide Vittone}
\address{Davide Vittone: Scuola Normale Superiore\\ Piazza dei Cavalieri 7\\ 56126 Pisa - Italy} 
\email{d.vittone@sns.it}
\thanks{D.V. is supported by MIUR, Italy, GNAMPA of INDAM and Scuola Normale Superiore, Italy. Part of the work was done while D.V. was a visitor at the University of Trento. He wishes to thank the Department of Mathematics for its hospitality.}

\date{\today}
\maketitle


\section{Introduction}
In this paper we deal with some problems concerning minimal hypersurfaces in Carnot-Carath\'eodory (CC) structures. More precisely we will introduce a general {\it calibration method} in this setting and we will study the {\it Bernstein problem} for entire regular intrinsic minimal graphs in a meaningful and simpler class of CC spaces, i.e. the Heisenberg group $\Hn$. In particular we will positevily answer to the Berstein problem in the case $n=1$ and we will provide counterexamples when $n\ge 5$.

Here hypersurface simply means topological codimension 1 surface. The notion of intrinsic graph has been recently introduced in \cite{FSSC04} in the setting of a Carnot group and deeply studied in the case of hypersurfaces in \cite{ASCV}  even if it was already implicitly exploited in \cite{FSSC}.

The calibration method presented here is a refinement of an unpublished result due to L. Ambrosio in \cite{ambrosio}. The Bernstein problem in the framework of CC spaces (see Definition \eqref{definizione d} below) has been recently studied for the first Heisenberg group $\H^1$  with suitable assumptions we will precise below (see \cite{GP1},\cite{CHMY}, \cite{CHw},\cite{RR2} and \cite{DGN2006}).

Many other topics of Geometric Measure Theory (GMT) have been studied in the setting of these structures (see, for instance, \cite{CDG}, \cite{Fragalwhe},\cite{FSSC1},\cite{GN},\cite{FSSC},\cite{ambrosio3}, \cite{ambrosio2},\cite{FSSC step2} and \cite{magn3} where an interesting survey of this argument can be found). In particular, given a family $X=(X_1,\dots,X_m)$ of Lipschitz continuous vector field on $\R^n$, it is possible to define the notion of $X$-perimeter $\norma\partial E\norma_X(\Omega)$ of a measurable $E\subset\R^n$ in open set $\Omega\subset\R^n$ (see Definition \eqref{defXperim}). Several interesting properties of sets of finite $X$-perimeter have been investigated (see \cite{CDG}, \cite{birmos}, \cite{FSSC1}, \cite{GN}) and in particular {\it{isoperimetric type inequalities}} have been proved for the $X$-pe\-ri\-me\-ter $\norma\partial E\norma_X$ associated to a measurable set $E\subset\Rn$ of {\it{locally finite $X$-perimeter}}, that is
\begin{equation*}
\norma\partial E\norma_X(\Omega)<\infty
\label{perimetro}
\end{equation*}
for every open bounded set $\Omega\subset\Rn$ (see \cite{pansucras}, \cite{Fragalwhe} and \cite{GN}). In particular an interesting existence's result of minimizing sets of locally $X$-finite perimeter has been obtained in \cite{GN}, generalizing similar De Giorgi's results for the Euclidean perimeter measure $\norma\partial E\norma:=\norma\partial E\norma_X$ with $X=\nabla=(\partial_1,...,\partial_n)$ (see \cite{Giusti}).

Moreover, introducing a suitable intrinsic notions of rectifiability in a Carnot group of step 2, a counterpart of De Giorgi's result on the structure of sets of finite Euclidean perimeter (see, for instance, \cite{Giusti}) has been obtained in \cite{FSSC step2} (see also \cite{ambrosio3}, \cite{pauls},\cite{magn3},\cite{mattila2},).
 
In this paper we are going to study some properties for locally finite $X$-perimeter sets $E\subset\Rn$ which are locally minimizing in a fixed open subset $\Omega\subset\Rn$ of a CC-space $\Rn$, i.e.
\begin{equation}\label{minimizer}
\norma\partial E\norma_X(\Omega')\,\le\,\norma\partial F\norma_X(\Omega')\quad{\text{whenever}}\quad E\Delta F\Subset\Omega'
\end{equation}
for any  open subset $\Omega'\Subset\Omega$, for any measurable set $F\subset\Rn$. In the following we will simply call such a set  $E$ a {\it minimizer} for the $X$-perimeter measure in $\Omega$. 
 
Through a calibration technique we will first give, in Theorem \ref{calibrazioni}, a general simple method for the construction of minimizers for the $X$-perimeter measure. Then we will the study the Bernstein problem for minimizers $E$ in the whole $\Hn$; namely we will study the problem whether a set $E$ satisfying (\ref{minimizer}) with $\Omega=\Hn$ whose the boundary $S=\partial E$  is   a regular  entire ``graph'' has to be an ``half-space'' with respect to suitable notions of graph and half-space involving both algebraic and CC structure of $\Hn$ (see section 5).  

We refer to S.Berstein since he positively answered to this problem in the Euclidean setting, i.e. when $X=\nabla=(\partial_1,...,\partial_n)$ with $n=3$. More precisely it is well-known that if $\phi:\omega\to\R$ is a $\ci^2$ regular function on an open subset $\omega\subset\R^{n-1}$ and $E=\{(x',x_{n})=(x_1,\dots,x_{n-1},x_{n})\in\omega\times\R: x_{n}<\f(x')\}$ then 
\begin{equation}\label{euclidgrapharea}
\norma\partial E\norma_X(\omega\times\R)=\mathcal H^{n-1}(\partial E)=\int_{\omega}\sqrt{1+|\nabla \f|^2}d \mathcal{ L}^{n-1}
\end{equation}
where $\mathcal H^m$ and $\mathcal L^m$ respectively denote the $m$-dimensional Hausdorff measure and the $m$-dimensional Lebesgue measure in $\Rm$.

In particular whenever $E$ is a minimizer in $\Omega=\omega\times\R$ then $\f$ satisfies the classical {\it minimal surface equation}
\begin{equation}\label{EMSE}
\text{div}\left(\frac{\nabla \f}{\sqrt{1+|\nabla \f^2|}}\right)=0\quad\text{in}\quad\omega\,.
\end{equation}
When $n=3$ Bernstein proved that each $\f\in \ci^2$ entire solution of (\ref{EMSE}) (i.e. $\omega=\R^2$) has to be affine  (see Theorem \ref{bernsteincl} and \cite{Giusti}, chapter 17, for an interesting account of this celebrated problem in any dimension).

Let us recall both problems we will deal with have been already studied in the setting of first Heisenberg group $\H^1$ though by a different point of view and with different assumptions. Indeed the problem of minimizers was attacked in \cite{pauls1}, \cite{CHMY},\cite{HP}, \cite{RR2} and \cite{CHY} but only for sets $E$ with regular boundary and in the special setting of Heisenberg group $\Hn$. Our calibration technique, although very classic as inspiration (see, for instance, \cite{bugiahil} and \cite{Morgan} for its historical applications  to Calculus of Variations and GMT), extends previous results to general sets without any regularity assumptions and in the general framework of CC spaces (see section 2). 
 
The Bernstein problem in $\H^1$ was first studied in \cite{GP1} and \cite{CHMY} and in  the subsequent papers \cite{CHw}, \cite{GP2} and  \cite{RR2}. To introduce those and our results let us recall some definitions and the main properties of $\Hn$ concerning both its algebraic and CC structure (see \cite{capdanpautys} for a more complete introduction). For sake of simplicity in the introduction we will restrict ourselves to the first Heisenberg group $\H^1$.

Let us denote $\H^1=(\R^3,\cdot)$ with group law 
\begin{equation}\label{grouplaw}
(x,y,t)\cdot (x',y',t'):=(x+x',y+y',t+t'+2(yx'-xy'))
\end{equation}
and the family of non isotropic dilations
\begin{equation}\label{Hdilation}
\delta_{r}((x,y,t)):=(rx,ry,r^2t),  \text{ for }r>0.
\end{equation}
It is known that $\H^1$ is a Lie group of topological dimension $3$, its Lie algebra $\mathfrak{h}_1$ of the left invariant vector fields of $\H^1$ is (linearly) generated by
\begin{equation}\label{campiH1}
X_1=\frac{\partial}{\partial x}+2y\frac{\partial}{\partial t},\qquad X_2=\frac{\partial}{\partial y}-2x\frac{\partial}{\partial t }; \qquad T=\frac{\partial}{ \partial t},
\end{equation}
and the only non-trivial commutator relation is
\begin{equation*}
[X_1,X_2] = - 4T .
\end{equation*}
Let us denote in the following $X:=(X_1,X_2)$ and let $d_c$ be the CC distance defined in (\ref{definizione d}): it is well-known that $d_c$ is left-invariant with respect to the group law and homogeneous with respect to the dilations $(\delta_r)_r$. In what follows, we will write $\norma\partial E\norma_\heis$ rather than $\norma\partial E\norma_X$, and we will call it $\heis$-perimeter instead of $X$-perimeter.

Actually $\H^1$ can be endowed with an explicit simpler metric, $d_{\infty}$, equivalent to $d_c$ we will use through the paper. More precisely let us define the homogeneous norm in $\H^1$
\begin{equation}\label{Hnorm}
\norm {P}:=\max\left\{\sqrt {x^2 + y^2}, |t|^{1/2}\right\}\quad\text{if}\quad P=(x,y,t)
\end{equation}
then the metric $d_{\infty}$ we shall deal with will be defined as
\begin{equation}\label{distanced}
d_{\infty}(P,Q):=\norm{P^{-1}\cdot Q}.
\end{equation}
Moreover let us stress that Hausdorff dimension of $(\H^1,d_{\infty})$ is $Q:=4$.

The two key points we want to stress now are about the possible intrinsic notions of plane and graph in $\H^1$ equipped with its Lie group and CC structure. A general and more complete discussion of these topics in Carnot groups can be found in \cite{FSSC04} and \cite{ASCV}. 

The notion of plane arises quite evident taking into account the fundamental Rademacher's type theorem due to P. Pansu in \cite{pansu} in the framework of Carnot groups. Indeed Pansu proved the intrinsic differential (i.e. an intrisic linear functional) $L:\H^1\to \R$, approximating with respect to the distance $d_{\infty}$ at a.e. point in a fixed open set $\Omega\subset\H^1$, a given Lipschitz continuous function $f:\Omega\subset (\H^1,d_{\infty})\to\R$ has to be the form
\begin{equation*}
L(x,y,t)=a\,x\,+b\,y\quad\forall \,(x,y,t)\in\H^1
\end{equation*}
for suitable constants $a$, $b\in\R$. Then it is natural by this argument to call {\it plane} in $\H^1$ a set $V\subset\H^1$ level set of $L$, i.e. $V=\{(x,y,t): ax+by=c\}$ form some $c\in\R$, i.e it is a so-called {\it vertical plane}. Let us note such a set $V=P_0\cdot V_0$ for some $P_0\in V$ and $V_0:=\{(x,y,t)\in\H^1: ax+by=0\}$. Namely a plane can be seen as the left-transleted of a maximal subgroup $V_0$ of $\H^1$. This evidence has been later confirmed  by an argument of GMT too. In fact it was proved in \cite{FSSC} the existence of a tangent maximal subgroup of $\H^1$ in the sense of GMT for {\it $\H\,$-regular surfaces} (see section 2), a class of intrinsic regular surfaces in $\Hn$ which allowed to prove a counterpart of the classical De Giorgi's rectificability and divergence theorem in this setting (see also a generalization of this blow-up' s technique recently proposed in \cite{mattila2}).

Assuming this notion of vertical plane the Bernstein problem was studied in \cite{GP1} for $\ci^2$ regular sets
\begin{equation}\label{tgraph} 
E=\{(x,y,t)\in\R^3: t\,<\f(x,y)\}
\end{equation}
verifying (\ref{minimizer}) and with $S=\partial E$ without characteristic points. In particular a suitable {\it minimal surface equation} for $\f$ was obtained and its solutions have been called {\it H- minimal}. In particular it was proved that Bernstein problem fails for an entire $\ci^2$-regular H-minimal solutions $\phi$, i.e. there exists H-mimimal $\f:\R^2\to\R$ whose graph is not a vertical plane. On the other hand  it has been been recently proved in \cite{RR2} (see also \cite{CHY}) that for each $\ci^2$ regular  entire H-minimal solution $\f$ the set $E$ in (\ref{tgraph})) is a minimizer for the $\heis$-perimeter in the whole $\Omega=\H^1\equiv\R^3$.

In this paper we want to study the Bernstein probem in $\H^1$ assuming again the notion of vertical plane but replacing the notion of Euclidean graph with another one more intrinsic to the CC geometry which turned be very useful in GMT (see also \cite{DGN2006}). 

Let us recall that, according to \cite{ASCV}, two sets $S,\,E\subset\H^1$ are respectively called {\it $X_1$-graph} and  {\it $X_1$-subgraph} in $\H^1$ if  
\begin{equation}\label{Xgraph}
S=\left\{(0,\e,\t)\cdot \f(\e,\t)\,e_1:(\e,\t)\in\omega\right\}
\end{equation}
\begin{equation}\label{Xsubgraph}
E=\left\{(x,y,t)\in C_{X_1}(\omega): x<\f(y,t-2xy\right)\}
\end{equation}
for a suitable function $\f:\omega\subset\R^2\to\R$ where $C_{X_1}(\omega):=\{(x,y,t)\in\H^1: (y,t-2xy)\in\omega\}$ and $e_i\, (i=1,2,3)$ is the standard basis of $\R^3$. 

Let us point out that the notion of $X_1$-graph is not a pointless generalization. Indeed, for instance, there are $\H$-regular surfaces in $\H^1$ which are not Euclidean graphs (see \cite{FSSC04}, Example 3.9).

Moreover if $\f\in \ci^2(\omega)$ then $\partial E=S$ and the following area formula for $X_1$-graphs was proved in \cite{ASCV} (see Remark  2.23)
\begin{equation}\label{Hgrapharea}
\norma\partial E\norma_\heis(C_{X_1}(\omega))=\mathcal S_{\infty}^3(S)=c\,\int_\omega\sqrt{1+|W^\f\f|^2}\,d\mathcal L^2
\end{equation}
where $\Wf:=\dede{\e}-4\phi\dede\t$, $\mathcal S_{\infty}^3$ is the $3$-dimensional spherical Hausdorff measure induced on $\H^1$ by $d_{\infty}$ and $c$ a positive geometric constant.

By performing a simple first variation with respect to $\f$ of the functional in (\ref{Hgrapharea}) we can easily get (see section 3 and \cite{DGN2006}) the following {\it minimal surface equation for $X_1$-graphs} in $\H^1$
\begin{equation}\label{HMSE}
\Wf\left(\frac{\Wf\f}{\sqrt{1+|\Wf\f|^2}}\right)=0\qquad\text{in }\quad\omega\,.
\end{equation}
In \cite{GP2} interesting examples of (entire) $X_1$-graph satisfying (\ref{HMSE}) in the whole $\omega=\R^2$ which are not vertical planes have been provided. Therefore the only request to be an an entire solution of (\ref{HMSE}) does not yields a positive answer to the Bernstein problem.

On the other hand very recently in \cite{DGN2006} a remarkable example of a regular entire $X_1$-graph $S$ verifying (\ref{HMSE}) but not minimizing the $\heis$-perimeter, in the sense of (\ref{minimizer}) has been provided. Moreover it has been proposed the natural conjecture whether a $\ci^2$ regular set $E\subset\H^1$ as in (\ref{Xsubgraph}) with $\omega=\R^2$ and satisfying (\ref{minimizer}) with $\Omega=\H^1$ has to be a vertical half-space. 

In this paper, getting motivation and benifit of this interesting counterexample, we positively answer to this question (see Theorem \ref{bernsteinminim}) by means of a characterization of the entire solutions of (\ref{HMSE}) (see Corollary \ref{unicita}) and studying among them those satisfy the right sign of the second variation of functional in (\ref{Hgrapharea}) (see section 3.2).

Then we extend our study to the higher Heisenberg groups $\Hn$ with $n\ge 2$ (see section 5.2). In particular we will prove that when $n\ge 5$ there is a negative answer to Bernstein problem in $\Hn$ by constructing examples of regular minimizers  $E\subset\Hn\equiv\R^{2n+1}$ whose the boundary $S=\partial E$  is an entire intrinsic graph but not  a vertical hyperplane. Moreover we will study also the Bernstein problem in the cases n=2, 3 and 4 getting only some partial information but it  still remains open.

Finally  N. Garofalo informed us that similar topics have been studying in \cite{GS}.

\medskip

{\it Acknowledgements.} We are deeply grateful to L. Ambrosio who generously communicated us a preliminary version of Theorem \ref{calibrazioni} and for fruitful discussions on the subject. We also warmly thank S. Spagnolo who inspired us many results of section 4, P. Secchi for fruitful discussions on the arguments of the same section, R. Monti for some suggestions on minimal surfaces and N. Garofalo who kindly sent us a preliminary version of paper \cite{DGN2006} which gave us an important hint for the proof of Theorem \ref{bernsteinminim}. Eventually many thanks to B. Franchi and R. Serapioni for continuous discussions during the writing of this paper and to V. Chilla for some helpful advices.


\section{A general calibration method for the $X$-perimeter and applications}
Given a family $X=(X_1,...,X_m)$ of Lip\-schitz continuous vector fields $X_j(x)=\sum_{i=1}^nc_{ij}(x)\partial_i$ ($j=1,...,m$) with $c_{ij}\in$ Lip$(\Rn)$ ($j=1,...,m$, $i=1,...,n$), we call {\it{subunit}} a Lipschitz continuous curve $\gamma:[0,T]\longrightarrow\Rn$ such that
\begin{equation}
\dot{\gamma}(t)=\sum_{j=1}^m a_j(t)X_j(\gamma(t))\quad{\textrm{and}}\quad\sum_{j=1}^m a_j^2(t)\leq1 \quad{\textrm{ for a.e.}}t\in[0,T],\label{gamma punto}
\end{equation}
with $a_1,...,a_m$ measurable coefficients. Then define the CC distance between the points $x,y\in\Rn$ as
\begin{equation}
\begin{array}{ll}
d_c(x,y)=\inf\{T\geq0:&{\textrm{there exists a subunit path }}\gamma:[0,T]\rightarrow\Rn\\
&{\textrm{such that }} \gamma(0)=x{\textrm{ and }}\gamma(T)=y\}.
\end{array}
\label{definizione d}
\end{equation}
If the above set is empty put $d_c(x,y)=+\infty$. If $d_c$ is finite on $\Rn$, i.e. $d_c(x,y)<\infty$ for every $x,y\in\Rn$, it turns out to be a metric on $\Rn$ and the metric space $(\Rn,d)$ is called CC space (see, for instance, \cite{gromov2} and \cite{Montgom}). In particular we shall generally assume the following connectivity condition
\begin{equation}\label{cond(C)}
\text{$d_c$ is finite and the identity map ${\mathrm{Id}}:(\Rn,d_c)\to(\Rn,|\cdot|)$ is a homeomorphism.}
\end{equation}

Whenever $\Omega$ is an open subset of $\R^n$ and $f\in \ci^1(\Omega)$ we set 
$$Xf=(X_1 f,\dots,X_m f),$$
whereas if $\vf=(\vf_1,\dots,\vf_m)\in\ci^1_c(\Omega;\Rm)$ we put
\begin{equation}\label{defdivX}
\div_X\vf:=\sum_{j=1}^m X_j^\ast\vf_j,
\end{equation}
where $X_j^\ast$ is the adjoint operator of $X_j$ in $L^2(\Rn)$. Given a measurable subset $E\subset\R^n$ we define the $X$-perimeter of $E$ in $\Omega$ as
\begin{equation}\label{defXperim}
\norma\partial E\norma_X(\Omega):=\sup\left\{\int_E \div_X\vf:\vf\in\ci^1_c(\Omega,\Rm),|\vf|\leq 1\right\}\,;
\end{equation}
alternatively, we can define $\norma\partial E\norma_X(\Omega)$ as the total variation in $\Omega$ of $X\chi_E$.

It is well-known that by Riesz representation Theorem we can identify $\norma\partial E\norma_X$ as a Radon measure on $\Omega$ for which there exists an unique $\norma\partial E\norma_X$-measurable function $\nu_E:\Omega\to \Rm$ such that
$$\begin{array}{l}
|\nu_E|_{\Rm}=1\quad \norma\partial E\norma_X\text{-a.e. in }\Omega\vspace{.1cm}\\
\displaystyle\int_E\div_X\vf\,d\mathcal L^n= -\int_\Omega\langle\vf,\nu_E\rangle_{\Rm}\,d\norma\partial E\norma_X \quad\text{for all }\vf\in\ci^1_c(\Omega,\Rm).
\end{array}$$
In the following we will call $\nu_E$ {\it horizontal inward normal to} E (see \cite{FSSC1}). 

Finally, we will say that $E$ is a minimizer for the $X$-perimeter in $\Omega$ if
\begin{equation*}
\norma\partial E\norma_X(\Omega')\,\le\,\norma\partial F\norma_X(\Omega')
\end{equation*}
for any open set $\Omega'\Subset\Omega$ and any measurable $F\subset\Rn$ such that $E\Delta F\Subset\Omega'$.

The following result is a refinement of one due L. Ambrosio and it extends the classical calibration method giving sufficient conditions for a Borel set $E\subset\Rn$ to be minimizer of $X$-perimeter.

\begin{theo}\label{calibrazioni}
Let $\Omega\subset\Rn$ be an open set, let $X_1,\dots,X_m$ be a family of Lipschitz continuous vector fields in $\Omega$ and let $E$ be a set of locally finite $X$-perimeter in $\Omega$. 

Suppose there are two sequences $(\Omega_h)_h$ and $(\nu_h)_h$ such that
\begin{itemize}
\item[(i)]$\Omega_ h\subset\Omega$ is open, $\Omega_h \Subset\Omega_{h+1}$, $\cup_{h=1}^{\infty}\Omega_h=\Omega$;
\item[(ii)] $\nu_h\in \ci^1(\Omega;\Rm)$,  $|\nu_h(x)|_{\Rm}\leq 1$ for all $x\in\Omega$, $h$;
\item[(iii)] $\div_X \nu_h=0$ in $\Omega_h$ for each $h$;
\item[(iv)] $\nu_h(x)\to\nu_E(x)\ \norma\partial E\norma_X$-a.e. $x\in\Omega$.
\end{itemize}
Then $E$ is a minimizer for the $X$-perimeter in $\Omega$.
\end{theo}
\begin{proof}
Fix an open  set $\Omega'\Subset\Omega$ and a measurable set $F\subset\Rn$ such that $E\Delta F\Subset\Omega'$. Let $\Omega''$ be any open set such that $E\Delta F\Subset\Omega''\Subset\Omega'$. Let $\bar h$ and $\psi\in\ci^1_c(\Omega')$ be such that $\Omega'\subset\Omega_{\bar h}$, $0\leq\psi\leq 1$ and 
\begin{equation}\label{calibrazioni1}
\Omega''\Subset\{\psi=1\}\Subset \Omega'\Subset\Omega\,.
\end{equation}
Now notice that for each $h>\bar h$
\begin{equation}\label{calibrazioni2}
\int_\Omega\langle\nu_h\,\psi,\nu_E\rangle_{\Rm}\,d\norma\partial E\norma_X =\int_\Omega\langle\nu_h\,\psi,\nu_F\rangle_{\Rm}\,d\norma\partial F\norma_X
\end{equation}
Indeed by (\ref{calibrazioni1}) and $(iii)$
\begin{multline*}
\int_\Omega\langle\nu_h\,\psi,\nu_E\rangle_{\Rm}\,d\norma\partial E\norma_X -\int_\Omega\langle\nu_h\,\psi,\nu_F\rangle_{\Rm}\,d\norma\partial F\norma_X \\
= -\int_{\Omega' }(\chi_E-\chi_F)\div_X(\psi\nu_h)d\mathcal L^n = -\int_{\Omega''}(\chi_E-\chi_F)\div_X\nu_hd\mathcal L^n=0
\end{multline*}
By (\ref{calibrazioni2})
\begin{equation*}
\norma\partial F\norma_X(\Omega')\geq \left |\int_\Omega\langle\nu_h\,\psi,\nu_F\rangle_{\Rm}\,d\norma\partial F\norma_X\right | =\left |\int_\Omega\psi\langle \nu_h,\nu_E\rangle_{\Rm} d\norma\partial E\norma_X\right |\,.
\end{equation*}
By $(ii)$ and $(iv)$ and thanks to Lebesgue convergence theorem, as $h\to\infty$ we get
\begin{equation}\label{5.1.3}
\norma\partial F\norma_X(\Omega')\geq \int_{\Omega'}\psi\: d\norma\partial E\norma_X\ge \norma\partial E\norma_X(\Omega'').
\end{equation}

Let now increase $\Omega''\uparrow \Omega'$ in (\ref{5.1.3}) and we get the thesis.
\end{proof}



An interesting applications of Theorem \ref{calibrazioni} to non Euclidean regular sets $E$ can be given in a particular but significant class of CC spaces, namely the Carnot groups. Let us recall that a {\it Carnot group $\G$ of
step $k$} 
is a connected, simply connected Lie group whose Lie algebra ${\mathfrak{g}}$ admits a {\it step $k$ stratification}, i.e. there exist linear subspaces $\mathfrak g_1,...,\mathfrak g_k$ such that
\begin{equation}\label{stratificazione}
{\mathfrak{g}}=\mathfrak g_1\oplus...\oplus \mathfrak g_k, \quad [\mathfrak g_1,\mathfrak g_i]=\mathfrak g_{i+1},\quad
\mathfrak g_k\neq\{0\},\quad \mathfrak g_i=\{0\}{\textrm{ if }} i>k,
\end{equation}
where $[\mathfrak g_1,\mathfrak g_i]$ is the subspace of ${\mathfrak{g}}$ generated by the commutators $[X,Y]$ with $X\in \mathfrak g_1$ and $Y\in \mathfrak g_i$. Let $m_i=\dim(\mathfrak g_i)$, for $i=1,\dots,k$ and $h_i=m_1+\dots +m_i$ with $h_0=0$ and, clearly, $h_k=n$; we will also write $m$ instead of $m_1$. Choose a basis $e_1,\dots,e_n$ of $\mathfrak{g}$ adapted to the stratification, that is such that 
$$e_{h_{j-1}+1},\dots,e_{h_j}\;\text {is a base of}\; \mathfrak g_j\;\text{for each}\; j=1,\dots, k.$$

Let $X_1,\dots,X_{n}$ be the family of left invariant vector fields such that $X_i(e)=e_i$, where $e$ is the identity of the group. Given (\ref{stratificazione}), the subset $X=(X_1,\dots,X_m)$ generates by commutations all the other vector fields; we will refer to $X_1,\dots,X_m$ as {\it generating vector fields} of the group. The exponential map is a one to one map from $\mathfrak g$ onto $\G$, i.e. any $x\in\G$ can be written in a unique way as $x=\exp(x_1X_1+\dots+x_nX_n)$. Using these {\it exponential coordinates}, we identify $x$ with the n-tuple $(x_1,\dots,x_n)\in\R^n$ and we identify $\G$ with $(\Rn,\cdot)$ where the explicit expression of the group operation $\cdot$ is determined by the Campbell-Hausdorff formula (see \cite{folland}). In this coordinates $e=(0,\dots,0)$ and $(x_1,\dots,x_n)^{-1}=(-x_1,\dots,-x_n)$, and the adjoint operator in $L^2(\G)$ of $X_j$, $X_j^{*}$, $j=1,\dots,n$, turns out to be $-X_j$ (see, for instance,\cite{FSSC step2} Proposition 2.2).

Two important families of automorphism of $\G$ are the so called intrinsic translations and the intrinsic dilations of $\G$. For any $x\in\G$, the {\it (left) translation} $\tau_x:\G\to\G$ is defined as 
$$ z\mapsto\tau_x z:=x\cdot z. $$
For any $\lambda >0$, the {\it dilation} $\delta_\lambda:\G\to\G$, is defined as 
\begin{equation}\label{dilatazioni}
\delta_\lambda(x_1,...,x_n)=(\lambda^{\alpha_1}x_1,...,\lambda^{\alpha_n}x_n),
\end{equation} 
where $\alpha_i\in\N$ is called {\it homogeneity of the variable} $x_i$ in $\G$ (see \cite{FS} Chapter 1) and is defined as
\begin{equation}\label{omogeneita2}
\alpha_j=i \quad\text {whenever}\; h_{i-1}+1\leq j\leq h_{i},
\end{equation}
hence $1=\alpha_1=...=\alpha_m<\alpha_{m+1}=2\leq...\leq\alpha_n=k.$

The subbundle of the tangent bundle $T\G$ that is spanned by the family of vector fields $X=(X_1,\dots,X_m)$ plays a particularly important role in the theory, it is called the {\it horizontal bundle} $H\G$; the fibers of $H\G$ are 
$$ H\G_x=\mbox{span }\{X_1(x),\dots,X_m(x)\},\qquad x\in\G.$$
The sections of $H\G$ are called {\it horizontal sections}, a vector of $H\G_x$ is an {\it horizontal vector} while any vector in $T\G_x$ that is not horizontal is a vertical vector. Each horizontal section is identified by its canonical coordinates with respect to this moving frame $X_1(x),\dots,X_m(x)$. This way, an horizontal section $\vf$ is identified with a function $\vf =(\vf_1,\dots,\vf_m) :\rn{n} \rightarrow\R^m$, and this allows us to define $\div_X\vf$ as in \eqref{defdivX}.

Since the family $X$ satisfies condition \eqref{cond(C)}, it makes sense to define the Carnot-Carathéodory distance $d_c$ and the $X$-perimeter as in \eqref{definizione d} and \eqref{defXperim}. One can check that the CC distance $d_c$ is left invariant, i.e. $d_c(z\cdot x,z\cdot y)=d_c(x,y)$ for any $x,y,z\in\G$, and behave well with respect to the group dilations, i.e. $d_c(\d_\l x,\d_\l y)=\l d_c(x,y)$ for all $x,y\in\G$ and $\l>0$. The $n$-dimensional Lebesgue measure $\mathcal L^n$ is the Haar measure of the group, and the integer $Q:=m_1+2m_2+\dots+km_k$ is called homogeneous dimension of $\G$: it coincides also with the Hausdorff dimension of $\G$ with respect to $d_c$. Moreover, one has
$$\mathcal L^n(B(x,r))=r^Q \mathcal L^n(B(x,1))=r^Q \mathcal L^n(B(0,1)),$$
where $B(x,r)$ is the CC ball of center $x$ and radius $r$, and it can be checked that the $X$-perimeter measure is homogenous of degree $Q-1$, i.e.
$$\norma \partial(\d_\l E)\norma_X(\d_\l\Omega)=\l^{Q-1}\norma \partial E\norma_X(\Omega).$$

In Carnot groups one can refine Theorem \ref{calibrazioni} as follows:

\begin{theo}\label{calgruppi}
Let $\G=(\Rn,\cdot)$ be Carnot group. Let $E$, $\Omega$ be respectively a measurable and open set of $\Rn$, and denote by $\nu_E:\Omega\to\Rm$ the horizontal inward normal to $E$ in $\Omega$. Let us assume
\begin{itemize}
\item[(i)] $E$ has locally finite $X$-perimeter in $\Omega$;
\item[(ii)] $\div_X(\nu_E)=0\quad\text{in}\quad\Omega$ in distributional sense;
\item[(iii)] there exists an open set $\tilde\Omega\subset\Omega$ such that  $\norma\partial E\norma_X(\Omega\setminus\tilde\Omega)=0$ and $\nu_E\in \ci^0(\tilde\Omega)$.
\end{itemize}

Then $E$ is a minimizer of the $X$-perimeter in $\Omega$.
\end{theo}

Before the proof of Theorem \ref{calgruppi}, let us recall the classique technique of intrinsic convolution in homogeneous groups (see \cite{FS}), Proposition 1.20):
 
\begin{lemma}\label{lemmacal}
Let $\G=(\Rn,\cdot)$ be a Carnot group and let $\rho\in \ci^{\infty}(\Rn )$ be such that $0\le\rho\le 1$, $\int_{\Rn}\rho \,d\mathcal L^n=1$, $\spt\rho\subset B(0,1)$ and $\rho(x^{-1})=\rho(-x)=\rho(x)$ for all $x\in\Rn$. Let us denote 
\begin{equation}\label{mollifier}
\rho_{\eps}(x):=\eps^{-Q}\rho\left(\delta_{1/\eps}(x)\right),\qquad x\in\Rn\,;
\end{equation}
 
\begin{equation}\label{convolmollifier}
(\rho_{\eps}\star f)(x):=\int_{\Rn} \rho_{\eps}(y)\,f(y^{-1}\cdot\,x)d\mathcal L^n(y)=\int_{\Rn} \rho_{\eps}(x\cdot\, y^{-1})\,f(y)d\mathcal L^n(y)\,.
\end{equation}
Then
\begin{itemize}
\item[(i)] if $f\in L^p(\Rn),\ 1\le p<\infty$, then $\rho_{\eps}\star f\in \ci^{\infty}(\Rn)$ and $\rho_{\eps}\star f\to f$ in $L^p(\Rn)$ as $\eps\to 0$;
\item[(ii)] $\spt\ \rho_{\eps}\star f\subset B(0,\eps)\cdot\,\spt\ f$;
\item[(iii)] $X_j(\rho_{\eps}\star \vf)=\rho_{\eps}\star X_j\vf$ for any $\vf\in \ci_c^{\infty}(\Rn)$ and each $j=1,\dots,m$;
\item[(iv)] $\int_{\Rn}(\rho_{\eps}\star f)\,g\,d\mathcal L^n= \int_{\Rn}(\rho_{\eps}\star g)\,f\,d\mathcal L^n$ for every $f\in L^{\infty}(\Rn)$, $g\in L^{1}(\Rn)$;
\item[(v)]  if $f\in L^{\infty}(\Rn)\cap \ci^0(\Omega)$ for a suitable open set $\Omega\subset\Rn$ then $\rho_{\eps}\star f\to f$ uniformly on compact subsets of $\Omega$ as $\eps\to 0$.
\end{itemize}
\end{lemma}

\begin{proof}[Proof of Theorem \ref{calgruppi}] 
Let $\rho_{\eps}$ be the family of mollifiers introduced in Lemma \ref{lemmacal} and denote $\bar\nu:\Rn\to\Rm$ as $\bar\nu\equiv\nu$ in $\Omega$, $\bar\nu\equiv 0$ in $\Rn\setminus\Omega$. Let us denote
$$ \nu_{\eps}(x):=(\rho_{\eps}\star \bar\nu)(x)= \bigl((\rho_{\eps}\star \bar\nu_1)(x),\dots,(\rho_{\eps}\star \bar\nu_m)(x)\bigr)\quad x\in\Rn\,. $$
Let us begin to prove that for a fixed open set $\Omega'\Subset\Omega$
\begin{equation}\label{calgruppi1}
\int_{\Omega}\vf\,\div_X\nu_{\eps}\,d\mathcal L^n=0
\end{equation}
for every $\vf\in \ci_c^{\infty}(\Omega')$, $0<\eps<\frac{dist(\Omega',\Rn\setminus\Omega)}{2}$.
 
Fix $\vf\in \ci_c^{\infty}(\Omega')$; since $\vf_{\eps}:=\rho_{\eps}\star \vf \in \ci_c^{\infty}(\Omega)$ and $X_j$ ($j=1,\dots,m$) are self-adjoint, by Lemma \ref{lemmacal} $(i)$-$(iv)$ we can integrate by parts getting
$$\int_{\Omega}\div_X(\nu_{\eps})\,\vf\,d\mathcal L^n=-\int_{\Omega}\sum_{j=1}^m <\nu,X_j\vf_{\eps}>_{\Rm}\,d\mathcal L^n=0\,.$$
From (\ref{calgruppi1}) we get
\begin{equation}\label{calgruppi2}
\div_X(\nu_{\eps})=0\qquad\text{in }\Omega'
\end{equation}
for every open set $\Omega'\Subset\Omega$ provided $0<\eps<\frac{dist(\Omega',\Rn\setminus\Omega)}{2}$. 
 
Let $(\Omega_h)_h$ be a sequence of open subsets of $\Omega$ verifying $(i)$ of Theorem \ref{calibrazioni}. Then by (\ref{calgruppi2}) with $\Omega'\equiv\Omega_h$ there exists a sequence $\eps_h\to 0$ for which the sequence $\nu_h\equiv\nu_{\eps_h}$ satisfies the assumptions of Theorem \ref{calibrazioni}: indeed $(i)$-$(iii)$ therein are immediately satisfied. Let us prove $(iv)$.
 
By $(iii)$ and Lemma \ref{lemmacal} $(v)$ we get that $\nu_h\to\nu$ uniformly on compact subsets of $\tilde\Omega$ and then $(iv)$ of Theorem \ref{calibrazioni} follows.
\end{proof}

The simplest example of Carnot group is the Heisenberg group $\heis^n$: it is the Carnot group whose $(2n+1)$-dimensional Lie algebra $\mathfrak h_n$ is linearly generated by the fields $X_1,\dots,X_n,Y_1,\dots,X_n,T$, where the only non vanishing commutation relations are
$$[X_j,Y_j]=-4T\qquad\text{for all }j=1,\dots,n.$$
It follows that the algebra admits the stratification $\mathfrak h_n=\mathfrak g_1\oplus \mathfrak g_2$, where
$$\mathfrak g_1=\text{span }\{X_1,\dots,X_n,Y_1,\dots,X_n\}\quad\text{and}\quad\mathfrak g_2=\text{span }\{T\}.$$

Using exponential coordinates, throughout this paper we shall denote the points of $\heis^n$ by $P=[z,t]=[x+iy,t]$, $z\in\\ci^n$, $x,y\in\R^n$, $t\in\R$. If $P=[z,t]$, $Q=[\zeta,\tau]\in \heis^n$ and $r>0$, the group operation reads as
\begin{equation}\label{introgrouplaw}
P\cdot Q:=[z+\zeta,t+\tau+2\Im m(\langle z,\bar{\zeta}\rangle)]
\end{equation}
and the family of non isotropic dilations as
\begin{equation}\label{introHdilation}
\delta_{r}(P):=[rz,r^2t]. 
\end{equation}
The generators of the algebra can be written in coordinates as 
\begin{equation}\label{introdefcampi}
X_j =\frac{\partial}{\partial x_j}+2y_j\frac{\partial}{\partial t},\qquad Y_j=\frac{\partial}{\partial y_j}-2x_j\frac{\partial}{\partial t },\qquad \text{for}\;j=1,\dots,n;\qquad T=\frac{\partial}{ \partial t};
\end{equation}
we will frequently use the notation $X_j:=Y_{j-n}$ for $j=n+1,\dots,2n$. As we already agreed, we will define the Carnot-Carathèodory distance $d_c$ and the $X$-perimeter (which we will call $\heis$-perimeter) as the ones arising from the family $X=(X_1,\dots,X_{2n})$.

Moreover $\heis^n$ can be endowed with the homogeneous norm
\begin{equation}\label{introHnorm}
\|P\|_\infty:=\max\{|z|, |t|^{1/2}\}
\end{equation}
and the distance $d_{\infty}$ we shall deal with is defined as 
\begin{equation}\label{introdistanced}
d_{\infty}(P,Q):=\|P^{-1}\cdot Q\|_\infty.
\end{equation}
It is well-known that the distance $d_\infty$ is equivalent to $d_c$; it follows that the Hausdorff dimension of $(\heis^n,d_{\infty})$ is $Q=2n+2$, whereas its topological dimension is $2n+1$.

A real measurable function $f$ defined on an open set $\Omega\subset\heis^n$ is said to be of class $\ci^1_\heis(\Omega)$ if the distribution $\nabla_\heis f:=Xf=(X_1f,\dots,X_{2n}f)$ is represented by a continuous function. We shall say that $S\subset\heis^n$ is an $\heis$-regular hypersurface if for every $P\in S$ there exist a neighbourhood $U$ and a function $f\in\ci^1_\heis(U)$ such that $\nabla_\heis f\neq 0$ and $S\cap U=\{Q\in U:f(Q)=0\}$, and we will define the horizontal normal to $S$ at $P$ as
$$\nu_S(P):=-\frac{\nabla_\heis f (P)}{|\nabla_\heis f (P)|}.$$
Notice that a $\ci^1_\heis$ surface can be an highly irregular object from the viewpoint of Euclidean geometry, since it can be a fractal set (see \cite{kirsercas}). Since it is not restrictive, in the following we will deal only with surfaces $S$ which are locally zero level of functions $f\in\ci^1_\heis$ with $X_1f\neq0$.

A first important step in the study of $\heis$-regular surfaces is the Implicit Function Theorem \ref{DiniTheorem};  before stating it however we need to fix some notations. In the following we will identify the maximal subgroup $V_1=\{(x_1,\dots,x_{2n+1})\in\heis^n:x_1=0\}$ with $\R^{2n}$ through the map
\begin{eqnarray*}
\i &:& \R^2=\R_\e\times\R_\t\to V_1\\
&&(\e,\t)\longmapsto (0,\e,\t)
\end{eqnarray*}
if $n=1$, and 
\begin{eqnarray*}
\i &:& \R^{2n}=\R_\e\times\R_{v=(v_2,\dots,v_n,v_{n+2},\dots,v_{2n})}^{2n-2}\times\R_\t\to V_1\\
&&(\e,v,\t)\longmapsto (0,v_2,\dots,v_n,\e,v_{n+2},\dots,v_{2n},\t)
\end{eqnarray*}
for $n\geq 2$; we stress the strange choice made to denote the components of $v$, which however is justified by the structure of $\i$. We can then state (see \cite{FSSC}) the following 

\begin{theo}\label{DiniTheorem}{\bf[Implicit Function Theorem]}
Let $\Omega$ be an open set in $\Hn$, $0\in\Omega$, and let $f\in\ci^1_\heis(\Omega)$ be such that $X_1f(0) >0$, $f(0)=0$. Let
$$E:=\{P\in\Omega: f(P)<0\}\quad\text{and}\quad S:=\{P\in\Omega: f(P)= 0\};$$
then there exist $\de,h>0$ such that, if we put $I:=[-\delta,\delta]\times[-\delta,\delta]^{2n-2} \times[-\delta^2,\delta^2] \subset\R^{2n}_{\e,v,\t}$, $J:=\{(s,0,\dots,0)\in\heis^n:s\in [-h,h]\}$ and $\mathcal{U}:=\iota(I)\cdot J$, we have that
$$\begin{array}{l}
E \mbox{ has finite $\heis$-perimeter in } \mathcal U;\\
\partial E\cap\mathcal U=S\cap\mathcal U;\\ 
\nu_E(P)=\nu_S (P)\qquad\text {for all } P\in S\cap \mathcal U.
\end{array}$$

Moreover there exists a unique continuous function $\phi:I\rightarrow [-h,h]$ such that $S\cap\overline{\mathcal{U}}  = \F(I)$, where $\F$ is the map $I\ni A\mapsto\iota(A)\cdot(\f(A),0,\dots,0)\in\heis^n$, and the $\H$-perimeter has the integral representation
\begin{equation}\label{rappresintdiFSSC}
\norma\partial E\norma_\heis(\mathcal U) = \int_I \frac{|\nabla_\heis f|}{X_1f}(\F(A))\,d\mathcal L^{2n}(A).
\end{equation}
Finally the $\heis$-perimeter measure $\norma\partial E\norma_\heis$ coincides with $c(n)\mathcal S^{Q-1}_\infty\res S$, where $c(n)$ depends only on $n$ and $\mathcal S^{Q-1}_\infty$ is the $Q-1$ dimensional spherical Hausdorff measure arising from $d_\infty$.
\end{theo}

From now on, given a real continuous map $\f$ defined on an open set $\omega\subset\R^{2n}$, we will indicate with $\F$ the map $\omega\ni A\mapsto\iota(A)\cdot(\f(A),0,\dots,0)\in\heis^n$; therefore a natural question arise, i.e. to characterize all the $\f$ such that $S:=\F(\omega)$ is an $\heis$-regular surface. This problem has been solved in \cite{ASCV}: more precisely, if we set 
$$\begin{array}{l}
\displaystyle\Xt_j\f:=\frac{\partial\f}{\partial v_j}+2v_{j+n}\frac{\partial\f}{\partial\t},\quad \Yt_j\f:=\frac{\partial\f}{\partial v_{j+n}}-2v_j\frac{\partial\f}{\partial\t},\quad j=2,\dots n\vspace{.1cm}\\
\displaystyle\Yt_1\f:=\frac{\partial\f}{\partial\e}, \quad\Tt\f:=\frac{\partial\f}{\partial\t}\vspace{.1cm}\\
\displaystyle\Wf_{n+1}\f:=\frac{\partial\f}{\partial\e}-2\frac{\partial(\f^2)}{\partial\t}\vspace{.1cm}\\
\displaystyle\Wf\f:=\left\{\begin{array}{ll}
(\Xt_2\f,\dots,\Xt_n\f,\Wf_{n+1}\f,\Yt_{2}\f,\dots,\Yt_n\f) & \text{if }n\geq 2\\
\Wf_{2}\f & \text{if }n=1
\end{array}\right.
\end{array}$$
in distributional sense, then we have

\begin{theo}\label{teoASCV}
Let $\omega\subset\R^{2n}$ be an open set and let $\f:\omega\rightarrow\R$ be a continuous function. Let $\F:\omega\rightarrow\heis^n$ be the function defined by $\F(A):=\iota(A)\cdot(\f(A),0,\dots,0)$ and let $S:=\F(\omega)$. Then the following conditions are equivalent:
\begin{enumerate}
\item[(i)] $S$ is an $\heis$-regular surface and $\nu_S^{1}(P)<0$ for all $P\in S$, where $\nu_S(P)=(\nu_S^{1}(P),\dots,\nu_S^{2n}(P))$ is the horizontal normal to $S$ at a point $P\in S$;
\item[(ii)] the distribution $\Wf\f$ is represented by a continuous function and there exist a family $\{\f_\ep\}_{\ep>0}\subset\ci^1(\omega)$ such that, for any open set $\omega'\Subset \omega$, we have
\begin{equation}\label{equazione2delteo1}
\f_\ep\rightarrow\f\mbox{ and }W^{\f_\ep}\f_\ep\rightarrow \Wf\f\mbox{ uniformly on }\omega'.
\end{equation}
\end{enumerate}
Moreover, for all $P\in S$ we have 
\begin{equation}\label{formadellanormale}
\nu_S(P)=\left(-\frac{1}{\sqrt{1+|\Wf\f|^2}},\frac{\Wf\f}{\sqrt{1+|\Wf\f|^2}} \right) (\F^{-1}(P)),
\end{equation}
and
\begin{equation}\label{ASCVareaformula}
{\mathcal S}_{\infty}^{Q-1}(S)= c(n)\,\int_{\omega}\sqrt{1+|\Wf\f|^2}\,d{\mathcal L}^{2n}
\end{equation}
where ${\mathcal L}^{2n}$ denotes the Lebesgue measure on $\R^{2n}$.
\end{theo}

Throughout the paper, we will make a large use of Theorem \ref{teoASCV} and expecially of the area formula \eqref{ASCVareaformula}.

We have now all the tools to state some results about minimizers of the $X$-perimeter in CC spaces: for all of them our calibration results will be crucial.

\begin{ex}[Hypersurfaces with constant horizontal normal]\label{exconstnorm}
{\upshape Let $X=(X_1,\dots,X_m)$ a family of Lipschitz continuous vector fields on $\Rn$. Let $E\subset\Rn$ be a set of locally finite $X$-perimeter in an open set $\Omega\subset\Rn$ which admits a constant inward horizontal nornal $\nu_E$ in $\Omega$, i.e.  
$$\nu_E\equiv\nu_0\quad\norma\partial E\norma_X\text{-}\text{a.e. in }\,\Omega$$
for a suitable constant vector $\nu_0\in\Rm$. Then, thanks to Theorem \ref{calgruppi}, it is straightforward to check that $E$ is a minimizer for the $X$-perimeter.

Observe that many interesting questions, such as regularity and rectifiability, are open even in this quite simple class of sets: see e.g. Example \ref{exH1nonsmooth}.}
\end{ex}

\begin{ex}[$t$-graphs in $\heis^1$]
{\upshape Let $\G=\H^1\equiv\R^3$ and $X=(X_1,X_2)$ as in (\ref{campiH1}). Let $\f\in \ci^2(\omega)$ for a suitable open set $\omega\subset\R^2$ and let $E$ be as in (\ref{tgraph}). Let  $\Omega:=\omega\times\R$,  $S=\partial E\cap\Omega$ and let $C(S)=\{(x,y,t)\in\Omega:\,-\f_x(x,y)+2y=\f_y(x,y)+2x=0\}$ the set of so-called {\it characteristic points} of $S$. Then $C(S)$ is closed in $\Omega$ and it was proved in \cite{balogh} that $\mathcal H^2(C(S))=0$. On the other hand since $\mathcal S^3_{\infty}<<\mathcal H^2$ (see \cite{brsc}) and $\norma\partial E\norma_\heis<<S^3_{\infty}$ (see \cite{FSSC}) we get  
\begin{equation}\label{example2.1}
\norma\partial E\norma_\heis(\Omega\setminus\tilde\Omega)=0
\end{equation} 
where $\tilde\Omega:=\Omega\setminus C(S)$. 

A simple calculation shows the horizontal normal $\nu_E(x,y,t)$ is 
\begin{equation}\label{example2.2}
\nu_E (x,y,t)=-\frac{\nabla_\heis f(x,y,t)}{|\nabla_\heis f(x,y,t)|}=N(x,y)=(N_1(x,y),N_2(x,y))
\end{equation}
for each $(x,y,t)\in S\setminus C(S)$, where $f(x,y,t):=t-\f(x,y)$ if $(x,y,t)\in\Omega$ and
$$N(x,y):=\frac {(-\f_x(x,y)+2y,-\f_y(x,y)-2x)} {\sqrt{(-\f_x(x,y)+2y)^2+(\f_y(x,y)+2x)^2}}\quad(x,y,t)\in \tilde\Omega\,.$$
In this case the minimal surface equation has been studied in \cite{pauls1}, \cite{GP1} and \cite{CHMY} when $C(S)=\emptyset$ and it simply reads as
\begin{equation}\label{HMSEEG}
\div_\heis\nu_E=\div\ N=\frac{\partial}{\partial x}N_1+\frac{\partial}{\partial y}N_2=0 \quad\text{in}\quad\omega.
\end{equation}
In particular, whenever \eqref{HMSEEG} is satisfied in pointwise sense, we can apply at once Theorem \ref{calgruppi} obtaining that $E$ is a minimizer for the $\heis$-perimeter measure in $\Omega$.
 
Very recently the more delicate case when $C(S)\neq\emptyset$ has been studied in \cite{RR2} and \cite{CHY}. In particular in \cite{CHY} it has been proved that (\ref{HMSEEG}) holds in weak sense, i.e.
\begin{equation}\label{HMSEEGW}
\int_{\omega} <N,\nabla\vf>_{\R^2}\,d\mathcal L^2=0\quad\forall\vf\in \ci_c^1(\omega)\,,
\end{equation}
iff $\f$ is a minimizer of the area functional in $\H^1$ for Euclidean $t$-graph.

We can get a strong result by exploiting Theorem \ref{calgruppi}: in fact, if (\ref{HMSEEGW}) holds, by (\ref{example2.1}) and (\ref{example2.2}) we get that $E$ is a minimizer for $\heis$-perimeter in $\Omega$. In particular $E$ minimizes the $\H$- perimeter not only among sets whose boundary is an Euclidean $t$- graphs but in a very much larger class of competitors.

Eventually let us stress our technique applies to the case studied in \cite{RR2}, Theorem 5.3. Indeed in our setting $\omega=\R^2$, $\f(x,y)=2xy+ay+b$, 
$$
N(x,y)=\left (0,\frac{4x-a}{|4x-a|}\right )\quad(x,y,t)\in \tilde\Omega=\left \{(x,y,t):x\neq a/4\right\}
$$
being $a$, $b\in\R$ fixed constants. On the other hand by a simple calculation we get that (\ref{HMSEEGW}) holds and then $E$ is a minimizer of the $\H$- perimeter in $\Omega=\R^3$.
}
\end{ex}



\begin{ex}
{\upshape Let $\G=\H^1\equiv\R^3$, $X=(X_1,X_2)$ as in (\ref{campiH1}). Let $E$ be as in (\ref{Xsubgraph}) with $\f(\e,\t)=-\frac{\a\e\t}{1+2\a\e^2}$, $\omega=\R^2$ for given constants $\alpha>0$. It was proved in \cite{DGN2006} that $S=\partial E$ is an entire $X_1$-graph which is not minimizing for the $\heis$-perimeter measure in the whole $\Omega=\H^1\equiv\R^3$. Let us stress here the difference with Example 2.7: here in fact $S$ is not a minimizer for $\heis$-perimeter measure though it satisfies the intrinsic minimal surface equation (\ref{HMSE}) with $\omega=\R^2$.

On the other hand we can prove it is a minimizer in $\Omega=\R^3\setminus\{y=0\}$: indeed with a simple calculation we get
$$\nu_E(x,y,t):= -\frac{\nabla_\heis f(x,y,t)}{|\nabla_\heis f(x,y,t)|}= \frac{y}{|y|}\left(-\frac{y}{\sqrt{x^2+y^2}},\frac{x}{\sqrt{x^2+y^2}}\right),$$
where $f(x,y,t):=x+\alpha\,yt$. Moreover it easy to see that $\nu_E\in \ci^{\infty}(\Omega)$ and 
$$\div_\heis(\nu_E)=0\quad\text{in}\quad\Omega\,.$$
Therefore applying Theorem \ref{calgruppi} we get the thesis. It is not know whether $S$ is $\heis$-perimeter minimizing in a neighbourhood of any point $(0,y,0),y\in\R$.}
\end{ex}

\begin{ex}[Nonsmooth minimal surfaces in $\heis^1$]\label{exH1nonsmooth}
{\upshape We provide a way to product minimizers of the $\heis$-perimeter in $\heis^1$ whose regularity is not better than (Euclidean) Lipschitz.  
Examples with this regularity are also provided in \cite{CHY} for minimal Euclidean $t$- graphs and very recently S.  Pauls informed us of a work in progress on this subject.

Our  key idea is to construct a ``not too regular'' parametrization $\f:\omega\to\R$ such that $\bur\f=0$ on an open set $\omega\subset\R^2_{\e,\t}$: indeed this property ensures that the horizontal normal to the surface is constant $\nu\equiv X_1$, and we conclude by calibrating with a constant section $\nu\equiv X_1$. 

Since for a Lipschitz continuous $\f$ we have $\bur\f=(\partial_\e-4\f\partial_\t)\f$, the required condition is equivalent to $\f$ being constant along the integral curves of the vector field $\Wf$, i.e. to these integral curves being lines. Notice that, using the same notations of Section \ref{solint}, this is equivalent to look for (local) solutions of \eqref{DB} with initial conditions $A\equiv 0$.

We then start by fixing a Lipschitz function $\b:\R\to\R$, with $L:=$Lip $\b<+\infty$, which will give the ``initial value'' of $\f$ in the sense that we look for a $\f$ such that $\f(0,\cdot)=\b$ ($\b$ is simply the counterpart of the  function $B$ of Section \ref{solint}). Fix a point $(\e,\t)\in\R^2$, consider the integral curve of $\Wf$ passing through it and let $(0,t)$ be the point in which this line meets the $\t$-axis: the condition of $\f$ being constant along this line then becomes $-4\f(\e,\t)=-4\b(t)=\frac{\t-t}{\e}$, i.e.
\begin{equation}\label{condiniz}
\t=t-4\e\b(t).
\end{equation}
Consider the Lipschitz continuous map 
\begin{eqnarray*}
F &:& \R^2_{x,t}\to\R^2_{\e,\t}\\
&& (x,t)\longmapsto (x,t-4x\b(t));
\end{eqnarray*}
$F$ plays the role of the $F$ of Section \ref{solint} and the variable $t$ the one of $c$. Since $t\mapsto F(x,t)$ is coercive when $x<1/4L$, one can easily see that $F$ is invertible when $\e<1/4L$. If we put $F^{-1}(\e,\t)=(\e,t(\e,\t))$ it turns out that condition \eqref{condiniz} is equivalent to define $\f(\e,\t):=\b(t(\e,\t))$, where from now on we suppose 
$$(\e,\t)\in\omega:=\ \left]-\frac{1}{4L}, \frac{1}{4L}\right[\times\R;$$
observe that $\f$ has the same regularity of $\b$ (but no more since $\f(0,\t)=\b(\t)$).  

Finally, let us verify that $\bur\f\equiv 0$: as
$$\nabla F(x,t)=\left(\begin{array}{cc}
1 & 0\\
-4\b(t) & 1-4x\b'(t)
\end{array}\right)$$
holds almost everywhere, one must have 
$$\nabla F^{-1}(\e,\t)=\bigl(\nabla F(F^{-1}(\e,\t))\bigr)^{-1}=\left(\begin{array}{cc}
1 & 0\vspace{0.3cm}\\
\displaystyle\frac{4\b(t(\e,\t))}{1-4\e\b'(t(\e,\t))} & \displaystyle\frac{1}{1-4\e\b'(t(\e,\t))}
\end{array}\right)$$
a.e., and so
\begin{eqnarray*}
&& \bur\f(\e,\t)=\bigl(\partial_\e-4\b(t(\e,\t))\partial_\t\bigr)\b(t(\e,\t))\\
&=& \b'(t(\e,\t))\frac{\partial t(\e,\t)}{\partial\e}- 4\b(t(\e,\t))\b'(t(\e,\t))\frac{\partial t(\e,\t)}{\partial\t}\\ 
&=& \b'(t(\e,\t))\frac{4\b(t(\e,\t))}{1-4\e\b(t(\e,\t))} -\frac{4\b(t(\e,\t))\b'(t(\e,\t))}{1-4\e\b'(t(\e,\t))}=0
\end{eqnarray*}
holds in the sense of distribution.

We stress that all the maps $\f:\omega\to\R$ arising from the previous discussion effectively parameterize a $\ci_\heis^1$ surface; in fact (see Theorem \ref{teoASCV}) it is sufficient to find $\ci^\infty$ functions $\f_\ep:\omega\to\R$ such that
$$\begin{array}{l}
\f_\ep\to\f \qquad\text{locally uniformly on }\omega\\
W^{\f_\ep}\f_\ep\to 0\qquad\text{locally uniformly on }\omega
\end{array}$$
as $\ep\to 0$. 

Fix then (e.g. mollifying $\b$) a sequence $\b_\ep\in\ci^\infty$ such that Lip $\b_\ep\leq L$ and $\b_\ep\to\b$ locally uniformly in $\R$, and consider the maps $\f_\ep$ arising from the previous discussion but considering $\b_\ep$ instead of $\b$. By construction we have $W^{\f_\ep}\f_\ep\equiv 0$; moreover, $\f_\ep$ are well defined on all $\omega$ (since Lip $\b_\ep\leq L$) and it is not difficult to check that they converge locally uniformly to $\f$.

Observe that if $\b$ is not $\ci^1$, then the surface parameterized by $\f$ cannot be of class $\ci^1$, since its intersection with the plane $\{y=0\}$ is the line $\{(\b(t),0,t):t\in\R\}$ which is not $\ci^1$.

For instance, let us put $\b(t)=|t|$: it is not difficult to compute that the asso\-cia\-ted parametrization is
\begin{eqnarray*}
\f &:& ]-1/4,1/4[\:\times\:\R\to\R\\
&& (\e,\t)\longmapsto\left\{\begin{array}{ll}
\displaystyle\frac{\t}{1-4\e} & \mbox{if }\t\geq 0\vspace{0.2cm}\\
\displaystyle-\frac{\t}{1+4\e} & \mbox{if }\t<0.
\end{array}\right.
\end{eqnarray*}
The surface parameterized by this $\f$ is then a perimeter minimizing set of class $\ci_\heis^1$ but not $\ci^1$.}
\end{ex}


\section{First and second variation of the area functional for intrinsic graphs}
In this section we want to obtain first and second variation formulas of the area functional for intrinsic graphs; similar formulas have been obtained in \cite{DGN1} for general surfaces. We will study in Section \ref{solint} the structure of all entire stationary points (i.e. those functions with vanishing first variation), while a proper second variation formula (cfr. \eqref{varsec0.4}) will be central in the study of the Bernstein problem in $\heis^1$ (see Section \ref{bernstein1}).

\subsection{First variation of the area}
Let us fix a $\ci^1$ map $\f:\omega\to\R$, where $\omega$ is an open subset of $\R^{2n}$, and put
\begin{equation}\label{Efi}
E_\f:=\{\i(A)\cdot(s,0,\dots,0)\in\heis^n:A\in\omega\text{ and }s<\f(A)\}\subset C_{X_1}(\omega)
\end{equation}
where we $C_{X_1}(\omega)$ is the cylinder of base $\i(\omega)$ along $X_1$ defined by
$$C_{X_1}(\omega):=\i(\omega)\cdot\{(s,0,\dots,0)\in\heis^n:s\in\R\};$$
observe that $C_{X_1}(\omega)$ is an open neighbourhood of $S:=\F(\omega)$, where as usual $\F$ is the map $A\mapsto\i(A)\cdot(\f(A),0,\dots,0)$. 

Let us assume that $E_\f$ is a minimizer for the $\heis$-perimeter in $C_{X_1}(\omega)$, fix $\p\in\ci^\infty_c(\omega)$ and set $\f_s:=\f+s\p$; we can therefore consider the class of competitors $E_{\f_s}$, which are defined as in \eqref{Efi} (observe that $E\Delta E_{\f_s}\Subset C_{X_1}(\omega)$), and set
\begin{equation}\label{g(s)}
g(s):=\norma\partial E_{\f_s}\norma_\heis(C_{X_1}(\omega))=\int_\omega\sqrt{1+|\Wfs\f_s|^2}d\mathcal L^{2n}.
\end{equation}
The fact that $g(s)\geq g(0)$ for all $s\in\R$ implies that $g'(0)=0$. We recall that the operator $\Wf$ is given by
$$\begin{array}{ll}
\Wf=(\Xt_2,\dots,\Xt_n,\Yt_1-4\f\Tt,\Yt_2,\dots,\Yt_n) & \text{if }n\geq 2\\
\Wf=\Wf_2=\Yt_1-4\f\Tt & \text{if }n=1.
\end{array}$$
We will write $\Xt_j:=\Yt_{j-n}$ when $n+1\leq j\leq 2n$.

It is well known that
$$\Xt^\ast_j = - \Xt_j\text{ for all }2\leq j\leq 2n\qquad\text{and}\qquad\Tt^\ast=-\Tt$$
and one can check that
$$(\Wf_{n+1})^\ast\p = -\Wf_{n+1}\p+4\p\Tt\f\qquad\text{for all }\p\in\ci^\infty.$$
Therefore we have
\begin{eqnarray*}
\Wfs_{n+1}\f_s &=& \Yt_1\f + s\Yt_1\p - 4(\f+s\p)(\Tt\f+s\Tt\p)\\
&=& \Wf_{n+1}\f - s\:{(\Wf_{n+1})}^\ast\p - 4s^2\p\Tt\p
\end{eqnarray*}
and so
\begin{equation}\label{g(s)espl}
g(s)=\int_\omega\Biggl[ 1+\sum_{\substack{j=2\\ j\neq n+1}}^{2n}(\Xt_j\f+s\Xt_j\p)^2+\bigl(\Wf_{n+1}\f - s\:{(\Wf_{n+1})}^\ast\p - 4s^2\p\Tt\p\bigr)^2 \Biggr]^{1/2}d\mathcal L^{2n}.
\end{equation}
From now on we will write just $\sum_j$ to mean the sum on indices $j=2,\dots,2n, j\neq n+1$; when $n=1$ the previous formula and the following ones are to be understood by ``erasing'' all sums of this type.

Starting from \eqref{g(s)espl} it is not difficult to compute
\begin{equation}\label{varprima}
g'(s)=\int_{\omega} \frac{\sum_j \Xt_j\f_s\:\Xt_j\p + \Wfs_{n+1}\f_s\:(-{(\Wf_{n+1})}^\ast\p-8s\p\Tt\p)}{\sqrt{1+|\Wfs\f_s|^2}}d\mathcal L^{2n}
\end{equation}
and in particular
\begin{equation}\label{varprima0}
g'(0)=\int_{\omega} \frac{\sum_j \Xt_j\f\:\Xt_j\p - \Wf_{n+1}\f\:{(\Wf_{n+1})}^\ast\p}{\sqrt{1+|\Wf\f|^2}}d\mathcal L^{2n}
\end{equation}
The Euler equation for stationary points of the area functional is then 
\begin{equation}\label{ESM3}
\Wf\cdot\frac{\Wf\f}{\sqrt{1+|\Wf\f|^2}} =0\qquad\text{on }\omega,
\end{equation}
where the previous equality must be understood in distributional sense.

\subsection{Second variation of the area}
If $\f\in\ci^1$ from \eqref{varprima} we can compute
\begin{eqnarray}
g''(s) &=& \int_\omega\frac{1}{1+|\Wfs\f_s|^2} \Biggl\{ \sqrt{1+|\Wfs\f_s|^2}\times \nonumber\\
&& \times\Biggl[\sum_j (\Xt_j\p)^2+\bigl({(\Wf_{n+1})}^\ast\p+8s\p\Tt\p\bigr)^2 - 8\p\Tt\p\Wfs_{n+1}\f_s \Biggr]+ \nonumber\\
&& -\Biggl[\frac{\Bigl[\sum_j \Xt_j\f_s\:\Xt_j\p + \Wfs_{n+1}\f_s\bigl(-{(\Wf_{n+1})}^\ast\p-8s\p\Tt\p\bigr)\Bigr]^2}{\sqrt{1+|\Wfs\f_s|^2}} \Biggr] \Biggl\}d\mathcal L^{2n}
\end{eqnarray}
and so
\begin{equation}\label{varsec0}
g''(0) = \int_\omega \frac{(1+|\Wf\f|^2)\left[ |{\Wf}^\ast\p|^2 - 8\p\Tt\p\Wf_{n+1}\f \right]- \Bigl(\Wf\f\cdot {\Wf}^\ast\p \Bigr)^2} {\left[1+|\Wf\f|^2\right]^{3/2}}d\mathcal L^{2n}
\end{equation}
where we put 
$$\begin{array}{ll}
{\Wf}^\ast\p:=\Bigl(\Xt_2^\ast\p,\dots,\Xt_n^\ast\p,{(\Wf_{n+1})}^\ast\p,\Xt_{n+2}^\ast\p,\dots,\Xt_{2n}^\ast\p\Bigr) & \text{if }n\geq 2\\
{\Wf}^\ast\p:={(\Wf_2)}^\ast\p & \text{if }n=1;
\end{array}$$
the fact that $E_\f$ is a minimizer implies that $g''(0)\geq 0$ for all $\p\in\ci^1_c(\omega)$.

Notice that when $n=1$ formula \eqref{varsec0} for the second variation reads as
\begin{equation}\label{varsec0.1}
g''(0) = \int_\omega \frac{|{\Wf}^\ast\p|^2 - 8\p\Tt\p\Wf\f(1+|\Wf\f|^2)} {\left[1+|\Wf\f|^2\right]^{3/2}}d\mathcal L^{2};
\end{equation}
in particular when $\Wf\f\equiv0$ one has $g''(0)\geq 0$ for all $\ci^1_c(\omega)$. 

If we suppose $\f\in\ci^2$ we can further exploit \eqref{varsec0.1} as
\begin{eqnarray}
g''(0) &=& \int_\omega \frac{|{\Wf}^\ast\p|^2 - 4\Tt(\p^2)\Wf\f(1+|\Wf\f|^2)}{\left[1+|\Wf\f|^2\right]^{3/2}}d\mathcal L^{2}\nonumber\\
&=& \int_\omega\Biggl[ \frac{|{\Wf}^\ast\p|^2}{\left[1+|\Wf\f|^2\right]^{3/2}} + 4\p^2\Tt\left(\frac{\Wf\f}{\left[1+|\Wf\f|^2\right]^{1/2}} \right)\Biggr]d\mathcal L^{2}.\label{varsec0.2}
\end{eqnarray}
We will see in Section \ref{solint} that if $n=1$ and $\f$ is a stationary point of the area functional, i.e. if $\f$ solves \eqref{ESM3}, then 
\begin{equation}\label{Wlapl0}
(\Wf)^2\f=0
\end{equation}
and thanks to this fact integrating by parts the first term of \eqref{varsec0.2} becomes
\begin{equation}\label{note(3)}
\int_\omega \frac{|{\Wf}^\ast\p|^2\,d\mathcal L^{2}}{\left[1+|\Wf\f|^2\right]^{3/2}} = \int_\omega \p\Wf\left( \frac{{\Wf}^\ast\p}{\left[1+|\Wf\f|^2\right]^{3/2}} \right) d\mathcal L^{2}= \int_\omega \p\frac{\Wf{\Wf}^\ast\p}{\left[1+|\Wf\f|^2\right]^{3/2}}d\mathcal L^{2}.
\end{equation}
Since
\begin{eqnarray*}
(\Wf{\Wf}^\ast-{\Wf}^\ast\Wf)\p &=& \Wf\bigl(-\Wf+4\Tt\f\text{Id}\bigr)\p- \bigl(-\Wf+4\Tt\f\text{Id}\bigr)\Wf\p\\
&=& 4\p\:\Wf\Tt\f
\end{eqnarray*}
we can rewrite \eqref{note(3)} as
\begin{eqnarray}
\int_\omega \frac{|{\Wf}^\ast\p|^2}{\left[1+|\Wf\f|^2\right]^{3/2}}d\mathcal L^{2} &=& \int_\omega\p\frac{{\Wf}^\ast\Wf\p + 4\p\:\Wf\Tt\f}{\left[1+|\Wf\f|^2\right]^{3/2}}d\mathcal L^{2}\nonumber\\
&=&\int_\omega \Biggl[\frac{(\Wf\p)^2}{\left[1+|\Wf\f|^2\right]^{3/2}} + 4\p^2\frac{\Wf\Tt\f}{\left[1+|\Wf\f|^2\right]^{3/2}}\Biggr] d\mathcal L^{2}\label{note(4)}
\end{eqnarray}
where we used \eqref{Wlapl0} again. Therefore \eqref{varsec0.2} becomes
\begin{eqnarray}
\vcenter{\begin{eqnarray*}
&&g''(0) = \int_\omega \Biggl\{\frac{(\Wf\p)^2}{\left[1+|\Wf\f|^2\right]^{3/2}} + 4\p^2\left[ \frac{\Wf\Tt\f}{\left[1+|\Wf\f|^2\right]^{3/2}} +  \Tt\left(\frac{\Wf\f}{\left[1+|\Wf\f|^2\right]^{1/2}}\right)\right] \Biggr\}d\mathcal L^{2}\nonumber\\
&&\hphantom{g''(0)}= \int_\omega \Biggl\{\frac{(\Wf\p)^2}{\left[1+|\Wf\f|^2\right]^{3/2}}+  \\
 &&\hphantom{g''(0)=}  +4\p^2\left[ \frac{\Wf\Tt\f}{\left[1+|\Wf\f|^2\right]^{3/2}} + \frac{[1+|\Wf\f|^2]\Tt\Wf\f-|\Wf\f|^2\Tt\Wf\p} {\left[1+|\Wf\f|^2\right]^{3/2}} \right]\Biggr\}d\mathcal L^{2}\nonumber\\
&&\hphantom{g''(0)} = \int_\omega \frac{(\Wf\p)^2 +4\p^2\:[\Wf\Tt\f+\Tt\Wf\f]}{\left[1+|\Wf\f|^2\right]^{3/2}}d\mathcal L^{2}
\nonumber\end{eqnarray*}}
\label{varsec0.3}\end{eqnarray}

Finally, we have
$$\Wf\Tt\f=\f_{\e\t}-4\f\f_{\t\t}=\Tt\Wf\f+4(\Tt\f)^2$$
and so from \eqref{varsec0.3} we can also write
\begin{equation}\label{varsec0.4}
g''(0) = \int_\omega \frac{(\Wf\p)^2 +8\p^2\:[\Tt\Wf\f+2(\Tt\f)^2]}{\left[1+|\Wf\f|^2\right]^{3/2}}d\mathcal L^{2}.
\end{equation}
Equation \eqref{varsec0.4} will be crucial in the proof of Theorem \ref{bernsteinminim}.


\section{Entire solutions of the minimal surface equation for intrinsic graphs in $\heis^1$}\label{solint}
In this section we will give a characterization (see Corollary \ref{unicita}) of all the entire $\ci^2$ solutions $\f:\R_{\e,\t}^2\to\R$ of the minimal surface equation for intrinsic graphs in $\heis^1$, i.e. of
\begin{equation}\label{MSE}
\Wf\left(\frac{\Wf\f}{\sqrt{1+|\Wf\f|^2}}\right)=0\qquad\text{in }\R^2;
\end{equation}
this result will be provide key tool to attack the Bernstein problem in $\heis^1$.

Observe that we can rewrite \eqref{MSE} also as
$$0 = \frac{(\Wf)^2\f\:\:\sqrt{1+|\Wf\f|^2}\:-\:\Wf\f\:\frac{\Wf\f\cdot(\Wf)^2\f}{\sqrt{1+|\Wf\f|^2}}}{1+|\Wf\f|^2}= \frac{(\Wf)^2\f}{(1+|\Wf\f|^2)^{3/2}}$$
which means that $\f$ is a solution of \eqref{MSE} if and only if it solves
\begin{equation}\label{MSE*}
(\Wf)^2\:\f=0\qquad\text{in }\R^2.
\end{equation}

Observe that \eqref{MSE*} is equivalent to a ``double'' Burgers'equation: in fact by performing the change of variables
\begin{eqnarray*}
G &:& \R^2_{x,t}\to\R^2_{\e,\t}\\
&& (x,t)\longmapsto (t,-4x),
\end{eqnarray*}
setting $u(x,t):=(\f\circ G)(x,t)=\f(t,-4x)$ and defining $L_u$ to be the operator
$$(L_u v)(x,t)=(v_t+uv_x)(x,t)\qquad(v\in\ci^1(\R^2)),$$
we get
$$(L_u(L_u u))(x,t)=((\Wf)^2\f)(t,-4x).$$
This means that we can restrict to consider the $\ci^2$ solutions $u$ of the ``double'' Burgers' equation
\begin{equation}\label{DB}
(L_u)^2u=0\qquad\text{in }\R^2
\end{equation}
(recall that $L_uu=0$ is the classical Burgers' equation). We will focus our attention on the problem \eqref{DB} rather than \eqref{MSE} or \eqref{MSE*}.

\subsection{Characteristic curves for entire solutions of \eqref{DB}}
Suppose $u$ is a $\ci^2$ solution of \eqref{DB} and let us consider the characteristic curves (see \cite{evans}) of the equation $L_uv=0$, i.e. for all $c\in\R$ the maximal solution $x=x(c,\cdot):I_c\to\R$ of the Cauchy problem
\begin{equation}\label{note1.1}
\left\{\begin{array}{l}
\dot x(c,t)=u(x(t,c),t)\\
x(c,0)=c.
\end{array}\right.
\end{equation}
Observe that from \eqref{DB} one gets $\frac{d}{dt}L_uu(x(c,t),t)=0$ and so
\begin{equation}\label{note1.2}
L_uu(x(c,t),t)=A(c)\qquad\text{for all }t\in I_c.
\end{equation}
Since
$$\frac{d}{dt}u(x(c,t),t)=\bigl(u_t(x(c,t),t)+u_x(x(c,t),t)\,\dot x(c,t)\bigr)=L_uu(x(c,t),t)=A(c)$$
we obtain
\begin{equation}\label{note1.3}
u(x(c,t),t)=A(c)t\,+\,B(c)\qquad\text{for all }t\in I_c,
\end{equation}
where we have set $B(c):=u(c,0)$. Equation \eqref{note1.3}, together with \eqref{note1.1}, gives
\begin{equation}\label{note1.4}
x(c,t)=\frac{A(c)}{2}t^2+B(c)t+c;
\end{equation}
in particular, $I_c=\R$.

We have therefore the following

\begin{theo}\label{teonote1.1}
Let $u$ be a $\ci^2$ solution of \eqref{DB} and for $c,t\in\R$ set $x(c,t):=\frac{A(c)}{2}t^2+B(c)t+c$, where $A(c)=L_uu(c,0)$ and $B(c)=u(c,0)$. Then for all $c,t$ we have
\begin{itemize}
\item[(i)] $u(x(c,t),t)=A(c)t+B(c)$;
\item[(ii)] $L_uu(x(c,t),t)=A(c)$;
\item[(iii)] $x(\cdot,t)$ is strictly increasing for any fixed time $t$;
\item[(iv)] for all $c\in\R$ we have one of $A'(c)=B'(c)=0$ or $B'(c)^2<2A'(c)$.
\end{itemize}
In particular, the family of characteristics $x(c,\cdot)$ are parabolas which cannot intersect themselves.
\end{theo}
\begin{proof}
We have already proved (i) and (ii); as for (iii), it will be sufficient to prove that, for every $t$,
\begin{equation}\label{note1.1.2}
x(c,t)\neq x(c',t)\qquad\text{if }c\neq c';
\end{equation}
in fact, were (iii) false, we could find $c<c'$ and $t'$ such that $x(c,t')\geq x(c',t')$, but since the characteristics are continuous and $x(c,0)=c<c'=x(c',0)$ we could find a $t$ between 0 and $t'$ such that \eqref{note1.1.2} does not hold.

Arguing by contradiction, let us assume that \eqref{note1.1.2} does not hold for some $c\neq c'$ and $t$; observe that, from (i) and (ii), one has
$$\begin{array}{l}
A(c)=L_uu(x(c,t),t)=A(c')\\
A(c)t+B(c)=u(x(c,t),t)=A(c')t+B(c')
\end{array}$$
whence $c=x(c,t)-\frac{A(c)}{2}t^2-B(c)t=c'$, which is a contradiction.

Notice that (iii) implies that
$$\frac{\partial x}{\partial c}(c,t)=\frac{A'(c)}{2}t^2+B'(c)t+1\geq 0$$
for all $c,t$, and this in turn implies $B'(c)^2\leq 2A'(c)$. Observe in particular that $A'(c)\geq 0$ and $\frac{\partial x}{\partial c}(c,t)\geq 0$.

In order to prove (iv), suppose by contradiction that for a certain $c$ we have $B'(c)^2= 2A'(c)\neq 0$. Let us differentiate (i) with respect to $c$ to get
$$\frac{\partial u}{\partial c}(x(c,t),t)=\frac{A'(c)t+B'(c)}{\frac{\partial x}{\partial c}(c,t)}= \frac{A'(c)\bigl(t+\frac{B'(c)}{A'(c)}\bigr)}{A'(c)\bigl(t+\frac{B'(c)}{A'(c)}\bigr)^2}=\frac{1}{t+\frac{B'(c)}{A'(c)}} \qquad\text{for all }t$$
which contradicts the hypothesis $u\in\ci^2(\R^2)$.
\end{proof}

\begin{obs}\label{remnote1.2}
{\upshape Observe that if $u$ is a $\ci^2$ solution of 
$$\left\{\begin{array}{l}
(L_u)^2u=0\\
u(x,0)=B(x)\\
L_uu(x,0)\equiv A\in\R
\end{array}\right.$$
then one must have also $B(x)\equiv B(0)=B$. In particular, Theorem \ref{teonote1.1} (i) implies that $u(x,t)=At+B$.}
\end{obs}

\begin{obs}\label{remnote1.3}
{\upshape Following the same proof of Theorem \ref{teonote1.1} (i), it is possible to prove that if $u$ is a $\ci^1$ solution of the Burgers' equation
$$L_uu=u_t+uu_x\equiv k$$
for a suitable constant $k\in\R$, then $B=u(\cdot,0)$ must be constant.}
\end{obs}

It is not difficult to extend the proof of Theorem \ref{teonote1.1} and get the following

\begin{theo}\label{teonote1.1bis}
Let $\Omega$ be an open set of $\R^2_{x,t}$ such that $\{(x,0):x\in\R\}\subset\Omega$, let $u\in\ci^2(\Omega)$ be a solution of \eqref{DB}, and let $A(c),B(c)$ and $x(c,t)$ be as in Theorem \ref{teonote1.1}. Suppose moreover that $\{(x(c,t),t):c,t\in\R\}\subset\Omega$. Then the statements (i)-(iv) of Theorem \ref{teonote1.1} still hold.
\end{theo}

From Theorem \ref{teonote1.1bis} we get the following uniqueness result for equation \eqref{DB} (see also \cite{CH}, Chap V, Section 7, and \cite{me}).

\begin{theo}\label{teonote1.4}
Let $u_0\in\ci^2(\R), u_1\in\ci^1(\R)$ be given functions and set $A:=u_0,B:=u_1+u_0u_0'$. Let $x(c,t):=A(c)t^2/2+B(c)t+c$ and suppose that 
\begin{equation}\label{UC}
\Omega=\{(x(c,t),t):c,t\in\R\}.
\end{equation}
Then there is at most one solution $u\in\ci^2(\Omega)$ of the problem
\begin{equation}\label{CDB}
\left\{\begin{array}{ll}
(L_u)^2u=0 & \text{in }\Omega\\
u(x,0)=u_0(x) & \forall\: x\in\R\\
u_t(x,0)=u_1(x) & \forall\: x\in\R.
\end{array}\right.
\end{equation}
\end{theo}
\begin{proof}
By Theorem \ref{teonote1.1bis} any solution $u\in\ci^2(\R^2)$ of \eqref{CDB} has to satisfy
$$u(x(c,t),t)=A(c)t+B(c);$$
however, hypothesis \eqref{UC} ensures that for all $(x,t)\in\Omega$ we can find a $c$ such that $x=x(c,t)$. This proves that $u$ is uniquely determined in $\Omega$ by $A$ and $B$, i.e. by $u_0$ and $u_1$.
\end{proof}

\begin{cor}\label{cornote1.4.1}
Let $u_0,u_1,A,B,x(t,c)$ and $\Omega$ be as in Theorem \ref{teonote1.4}, and suppose moreover that for all $c\in\R$ we have $A'(c)=B'(c)=0$ or $B'(c)^2<A'(c)$. Then
\begin{itemize}
\item[(i)] $\Omega$ is an open neigbourhood of the $x$-axis $\{(x,0):x\in\R\}$;
\item[(ii)] there is at most one solution $u\in\ci^2(\Omega)$ of the problem \eqref{CDB}.
\end{itemize}
\end{cor}
\begin{proof}
Observe that the map
\begin{eqnarray*}
F &:& \R^2\to\R^2\\
&& (c,t)\longmapsto (x(c,t),t)
\end{eqnarray*}
is regular and one-to-one; in particular, it is an open map and (i) follows. This means that condition \eqref{UC} of Theorem \ref{teonote1.4} is automatically fulfilled, and so (ii) must hold too.
\end{proof}


\begin{cor}\label{cornote1.5}
Under the same assumptions of Theorem \ref{teonote1.1} let us denote $l_1:=\lim_{c\to+\infty}A(c)$ (respectively $l_2:=\lim_{c\to-\infty}A(c)$). Then for any fixed $t\in\R$ we can conclude
\begin{equation}\label{tesicornote1.5}
\lim_{c\to+\infty} x(c,t)=+\infty\quad(\text{resp. }\lim_{c\to-\infty}x(c,t)=-\infty)
\end{equation}
if either $l_1\in\R$ (resp. $l_2\in\R$), or $l_1=+\infty$ (resp. $l_2=-\infty$) and one of the following conditions is satisfied:
\begin{eqnarray}
&& \liminf_{c\to+\infty} \frac{A(c)}{c}=0\quad \left(\text{resp. }\liminf_{c\to-\infty} \frac{A(c)}{c}=0\right) \label{note1.3.1}\\
&& \limsup_{c\to+\infty} \frac{A(c)}{c}=+\infty\quad \left(\text{resp. }\limsup_{c\to+\infty} \frac{A(c)}{c}=-\infty\right) \label{note1.3.2}\\
&& \liminf_{c\to+\infty} \left|\frac{B(c)}{\sqrt{cA(c)}}\right|<\sqrt{2}\quad \left(\text{resp. }\liminf_{c\to-\infty} \left|\frac{B(c)}{\sqrt{cA(c)}}\right|<\sqrt{2}\right).\label{note1.3.3}
\end{eqnarray}

In particular, when $\lim_{c\to+\infty} x(c,t)=+\infty$ and $\lim_{c\to-\infty}x(c,t)=-\infty$ we have that $x(\cdot,t):\R\to\R$ is an homeomorphism and $\Omega:=\{(x(c,t),t):c,t\in\R\}=\R^2$.
\end{cor}
\begin{proof}
Observe that for fixed $t\in\R$ and $c\neq 0$ one can write
\begin{equation}\label{note1.3.4}
x(c,t)=\sqrt{|c|}\left[ \frac{1}{2}\left(\frac{A(c)-A(0)}{\sqrt{|c|}}+\frac{A(0)}{\sqrt{|c|}}\right)t^2 + \left(\frac{B(c)-B(0)}{\sqrt{|c|}}+\frac{B(0)}{\sqrt{|c|}}\right)t +\sqrt{c} \right].
\end{equation}
Being $A$ increasing there exist
$$m_1:=\lim_{c\to+\infty} (A(c)-A(0))\qquad \bigl(\text{resp. }m_2:=\lim_{c\to-\infty} (A(c)-A(0))\bigr)$$
with $-\infty\leq m_2\leq 0 \leq m_1\leq + \infty$. Notice also that, using Theorem \ref{teonote1.1} (iv), one can get
\begin{equation}\label{note1.2.1}
|B(c)-B(0)|\leq\int_0^c|B'(s)|\,ds\leq\sqrt{2}\int_0^c \sqrt{A'(s)}\,ds\leq \sqrt{2\,|c|\,|A(c)-A(0)|}
\end{equation}
and this allows us to conclude when $l_1\in\R$ (resp. $l_2\in\R$), since in this case we have $m_1\in\R$ (resp. $m_2\in\R$) and so $x(c,t)\approx c$ for large (resp. small) $c$.

Instead, when $l_1=+\infty$, for large $c$ we can write
\begin{multline}\label{diseqschifosa}
x(c,t) = \sqrt{c(A(c)-A(0))}\left[ \frac{1}{2} \left(\sqrt{\frac{A(c)-A(0)}{c}}+\frac{A(0)}{\sqrt{c(A(c)-A(0))}}\right) t^2\right.+\\
\left. +\left( \frac{B(c)-B(0)}{\sqrt{c(A(c)-A(0))}} + \frac{B(0)}{\sqrt{c(A(c)-A(0))}}\right)t +\sqrt{\frac{c}{A(c)-A(0)}} \right]
\end{multline}
whence (using \eqref{note1.2.1} again) $\limsup_{c\to\infty} x(c,t)=+\infty$ in case \eqref{note1.3.1} or \eqref{note1.3.2} hold; however, this implies \eqref{tesicornote1.5} since $x(\cdot,t)$ is increasing. When $c\to-\infty$ we have instead
\begin{multline}\label{diseqschifosa2}
x(c,t) = \sqrt{c(A(c)-A(0))}\left[ -\frac{1}{2} \left(\sqrt{\frac{A(c)-A(0)}{c}}+\frac{A(0)}{\sqrt{c(A(c)-A(0))}}\right) t^2\right.+\\
\left. +\left( \frac{B(c)-B(0)}{\sqrt{c(A(c)-A(0))}} + \frac{B(0)}{\sqrt{c(A(c)-A(0))}}\right)t -\sqrt{\frac{c}{A(c)-A(0)}} \right]
\end{multline}
and we conclude analogously $\liminf_{c\to-\infty}x(c,t)=-\infty$, which is sufficient.

Instead if \eqref{note1.3.3} holds together with $l_1=+\infty$, we have a sequence $c_h\to+\infty$ such that
\begin{equation}\label{cgrande}
\frac{B'(c_h)^2}{2A'(c_h)}\leq (1-\ep)c_h\qquad\forall h;
\end{equation}
observe that the parabola $x(c_h,\cdot)$ reaches its minimum at $t=-\frac{B(c_h)}{A(c_h)}$ and so
$$x(c_h,t)\geq x(c_h,-\textstyle\frac{B(c_h)}{A(c_h)})=c_h-\displaystyle\frac{B'(c_h)^2}{2A'(c_h)}\geq\ep c_h \stackrel{h\to\infty}{\longrightarrow} +\infty$$
which, together with the fact that $x(\cdot,t)$ is increasing, proves \eqref{tesicornote1.5} when $c\to+\infty$. It is a little more complicated to get the thesis when $l_2=-\infty$ and $c\to-\infty$; however, as in \eqref{cgrande} we get a sequence $c_h\to-\infty$ such that
$$ - \frac{B'(c_h)^2}{2A'(c_h)} \leq \ep c_h -c_h\qquad\forall h$$
and so 
\begin{eqnarray*}
x(c_h,t) &=& \frac{A(c_h)}{2}\left( t+\frac{B(c_h)}{A(c_h)}\right)^2+\left(c_h-\frac{B'(c_h)^2}{2A'(c_h)}\right)\\
&\leq& \frac{A(c_h)}{2}\left( t+\frac{B(c_h)}{A(c_h)}\right)^2+\ep c_h
\end{eqnarray*}
which allows us to conclude since $A(c)\to-\infty$ as $c\to-\infty$.
\end{proof}

\begin{ex}\label{excaratt1}
{\upshape Let $A(c):=c/2$ and $B(c):=-c$; then it is easy to check that the family of characteristic curves for the related problem \eqref{DB} are
$$x(c,t)=(t-2)^2c/4.$$
Notice that $x(c,2)\equiv 0$, i.e. the thesis of Corollary \ref{cornote1.5} does not hold; here in fact \eqref{note1.3.3} is not fulfilled since
$$\lim_{c\to\pm\infty} \frac{B(c)}{\sqrt{cA(c)}}=\sqrt{2}.$$

Moreover, taking into account Theorem \ref{teonote1.1}, a global $\ci^2$ solution $u$ of \eqref{DB}, with $u(x,0)=-x$ and $L_uu(x,0)=x/2$, cannot exist. 
}
\end{ex}

\begin{ex}\label{excaratt2}
{\upshape Let $A(c)=c$ and $B(c)=\sqrt{2(1+c^2)}$, and let us consider the associated family of characteristic parabolas
$$x(c,t)=\frac{c}{2}t^2+\sqrt{2(1+c^2)}t+c.$$
Then for fixed $t$ we have
$$\frac{\partial x}{\partial c}(c,t)=\frac{t^2}{2}+\frac{\sqrt{2}c}{\sqrt{1+c^2}}t+1$$
which is (strictly) positive for any $c$: in particular, the family the characteristics cannot intersect, and in fact one has
$$B'(c)^2=2\frac{c^2}{1+c^2}<2=2A'(c).$$
Observe also that
$$\lim_{c\to\pm\infty}\left|\frac{B(c)}{\sqrt{cA(c)}}\right|=\sqrt{2}.$$

If we set $F(c,t):=(x(c,t),t)$ it is easy to see that the image $F(\R^2)$ is the open set
$$\Omega:=\R^2\setminus\left( \{(x,\sqrt{2}):x\leq 0 \}\cup\{(x,-\sqrt{2}):x\geq 0\} \right).$$
Indeed a simple calculation gives $F^{-1}(x,t)=(c(x,t),t)$ where 
$$c(x,t)=\left\{\begin{array}{l}
\displaystyle\frac{x(1+t^2/2)-\sqrt{2}|t|\sqrt{x^2+(1-t^2/2)^2}}{(1-t^2/2)^2} \quad \text{if }|t|\neq\sqrt{2}\\
\displaystyle\frac{x^2-4}{4x} \quad \text{if }t=\sqrt{2},x>0\text{ or } t=-\sqrt{2},x<0.
\end{array}\right.$$
We will see that $u(x,t):=A(c(x,t))t+B(c(x,t))$ is the unique solution of \eqref{DB} in $\Omega$ such that $L_uu(x,0)=A(x)$ and $u(x,0)=B(x)$.}
\end{ex}

\begin{ex}\label{excaratt3}
{\upshape If we require $B\equiv 0$, then the solution to \eqref{DB} with initial data $A,B$ is defined everywhere for any $\ci^2$ increasing function $A$. Obviously, even if it is possible to characterize it intrinsically as in Theorem \ref{teonote1.1} (i), in general it is not possible to give an explicit formula for the solution.}
\end{ex}

\subsection{Existence of entire solutions}
In the following theorem we provide an existence and uniqueness result for the equation \eqref{DB}.

\begin{theo}\label{teonote2.1}
Let $A,B\in\ci^2(\R)$ and for $c,t\in\R$ set
$$\begin{array}{l}
\displaystyle x(c,t):=\frac{A(c)}{2}t^2+B(c)t+c\vspace{0.1cm}\\
F:\R^2\ni(c,t)\longmapsto (x(c,t),t)\in\R\vspace{0.1cm}\\
\Omega:=F(\R^2)=\{(x(c,t),t):c,t\in\R\}
\end{array}$$
and suppose that
\begin{equation}\label{note2.1}
\text{for all } c\in\R\text{ one has either }A'(c)=B'(c)=0\text{ or }B'(c)^2<2A'(c).
\end{equation}
Then
\begin{itemize}
\item[(i)] $F$ is $\ci^2$ regular and one-to-one and, in particular, $\Omega$ is open;
\item[(ii)] if $F^{-1}(x,t):=(c(x,t),t),\ (x,t)\in\Omega$, then $u(x,t):=A(c(x,t))t + B(c(x,t))$ is the unique $\ci^2$ solution \eqref{DB} in $\Omega$ satisfying $L_uu(x,0)=A(x),\ u(x,0)=B(x)$.
\end{itemize}
\end{theo}
\begin{proof}
We begin by proving that the $\ci^2$ map $F:\R^2\to\R^2$ is one-to-one. By construction it is enough to prove that for any fixed $t$ the map $x(\cdot,t)$ is strictly increasing, and this is an easy consequence of \eqref{note2.1} which implies that
$$\frac{\partial x}{\partial c}(t,c)=\frac{A'(c)}{2}t^2+B'(c)+1$$
is strictly positive for any $c$. Being one-to-one and continuous, $F$ is also an open map, i.e. $\Omega\subset\R^2$ is open, and (i) is proved.

As for (ii), observe that the Jacobian matrix of $F$ is given by
$$JF(c,t)=\left(\begin{array}{cc}
\frac{A'(c)}{2}t^2+B'(c)+1 & A(c)t+B(c)\vspace{0.1cm}\\
0 & 1
\end{array}\right)$$
and so the Inverse Function Theorem implies that the Jacobian matrix of $F^{-1}$ is
\begin{multline*}
JF^{-1}(x,t) = \bigl( JF(F^{-1}(x,t))\bigr)^{-1}\\
= \frac{1}{\frac{A'(c(x,t))}{2}t^2+B'(c(x,t))+1}\left(\begin{array}{cc}
1 & -A(c(x,t))t-B(c(x,t))\vspace{0.1cm}\\
0 & \frac{A'(c(x,t))}{2}t^2+B'(c(x,t))+1
\end{array}\right).
\end{multline*}
Thus
\begin{equation}\label{note2.2}
\frac{\partial c}{\partial x}(x,t)=\frac{1}{\frac{A'(c(x,t))}{2}t^2+B'(c(x,t))+1}
\end{equation}
\begin{equation}\label{note2.3}
\frac{\partial c}{\partial t}(x,t)=-\frac{A(c(x,t))t+B(c(x,t))}{\frac{A'(c)}{2}t^2+B'(c)+1}
\end{equation}
and so one can compute
\begin{eqnarray*}
L_uu(x,t) &=& \bigl[A'(c(x,t))t + B'(c(x,t))\bigr]\frac{\partial c}{\partial t}(x,t) + A(c(x,t)) +\\
&& +\: \bigl[A'(c(x,t))t + B'(c(x,t))\bigr]\bigl[A(c(x,t))t + B(c(x,t))\bigr]\frac{\partial c}{\partial x}(x,t)\\
&=& A(c(x,t))
\end{eqnarray*}
and
$$L_u^2 u(x,t)= A'(c(x,t))\frac{\partial c}{\partial t}(x,t) + \bigl[A(c(x,t))t + B(c(x,t))\bigr]A'(c(x,t))\frac{\partial c}{\partial x}(x,t) =0.$$

Therefore $u$ is a solution of the given problem, and the proof is completed since uniqueness follows from Theorem \ref{teonote1.4}.
\end{proof}

\begin{cor}\label{unicita}
Suppose that $A,B\in\ci^2(\R)$ and that $u:\R^2\to\R$ is a $\ci^2$ entire solution of the problem 
$$\left\{\begin{array}{ll}
(L_u)^2u=0 \\
u(x,0)=B(x)\\
L_uu(x,0)=A(x)
\end{array}\right.$$
Let $\Omega,c(x,t)$ be as in Theorem \ref{teonote2.1}; then
$$u(x,t)=A(c(x,t))t+B(c(x,t))\qquad\text{for all }(x,t)\in\Omega$$
and $u$ is the unique solution in $\Omega$ of the same problem.
\end{cor}


\subsection{Examples of entire solutions of \eqref{DB}}

\begin{ex}\label{solint1}
{\upshape Let $A(c)=\a c$ ($\a>0$) and $B\equiv0$, then it is easy to see that in this case $\Omega=\R^2$; since $c(x,t)=\frac{2x}{2+\a t^2}$, the required solution of \eqref{DB} is given by
$$u(x,t)=\frac{2\a xt}{2+\a t^2}.$$
These solutions correspond to the maps $\f_{\a'}(\e,\t)=-\frac{\a'\e\t}{1+2\a'\e^2}$ (where $\a':=\a/4$) solutions of $(\Wf)^2\f=0$; it is not difficult to notice that the surfaces parameterized by $\f_{\a'}$ corresponds to $\{(x,y,t\in\heis^1:x=-\a'yt)\}$, which are deeply studied in \cite{DGN2006}: in particular (see Theorem 1.2 therein) it is proved that they are not $\heis$-perimeter minimizing (see also Theorem \ref{bernsteinminim}).}
\end{ex}

\begin{ex}\label{solint2}
{\upshape Let $B\equiv 0$ and choose a bounded, not constant and strictly increasing $A\in\ci^2$ ; then, if $\Omega$ and $c(x,t)$ are as in Theorem \ref{teonote2.1}, by Corollary \ref{cornote1.5} we have that $\Omega=\R^2$ and that $u(x,t):=A(c(x,t))t+B(c(x,t))$ is the unique entire solution of \eqref{DB}; moreover, $L_uu(x,t)=A(c(x,t))$ is bounded.


Observe that an analogous situation cannot occur in the Euclidean case: in fact (see \cite{Giusti}, Theorem 17.5), any smooth global solution $\p$ of the classical minimal surface equation with $\|\nabla \p\|_{L^\infty}<\infty$ must be linear. Here, instead, it happens that the map $\f$, which arises from the $u$ of this construction, solves \eqref{MSE}, is not linear (and, in particular, not of type \eqref{iperpiani}, see Section \ref{sezbernstein}) but is such that $\|\Wf\f\|_{L^\infty}<\infty $.}
\end{ex}


\section{The Bernstein problem in $\heis^n$}\label{sezbernstein}
The classical Bernstein problem is to find global functions $\p:\R^m\to\R$ solving the minimal surface equation
\begin{equation}\label{ESMC}
\text{div}\left(\frac{\nabla\p}{\sqrt{1+|\nabla\p|^2}}\right)=0
\end{equation}
and which are not affine functions, i.e. functions parameterizing hyperplanes or, which is the same, (translations of) maximal subgroups of $\R^{m+1}$. It is well known that this problem has been completely solved thanks to many contribution (see \cite{Giusti} for an interisting historical survey). Here we summarize these celebrated results in the following

\begin{theo}\label{bernsteincl}
Every $\p:\R^m\to\R$ which solves \eqref{ESMC} must be an affine function if $m\leq 7$; if $m\geq 8$ there are analytic functions $\p:\R^m\to\R$ solving \eqref{ESMC} which are not affine functions.
\end{theo}

\begin{obs}\label{ESMCmin}
{\upshape Notice that, through a standard calibration argument, one can prove that every Euclidean subgraph parameterized by an entire solution of \eqref{ESMC} is a minimizer for $X$- perimeter mesure in $\R^{n}$ provided $X=\nabla=(\partial_1,\dots,\partial_n)$ and $n=m+1$.}
\end{obs}

Let us recall the minimal surface equation for minimal $\heis$-graphs in $\heis^n$
\begin{equation}\label{ESM}
\Wf\cdot\left( \frac{\Wf\f}{\sqrt{1+|\Wf\f|^2}}\right)=0,
\end{equation}
where $\f:\R^{2n}\to\R$ is of class $\ci^2$. Observe that the ``affine'' functions given by
\begin{equation}\label{iperpiani}
\f(\e,v,\t)=c + \langle (\e,v),w \rangle_{\R^{2n-1}}
\end{equation}
for $c\in\R,w\in\R^{2n-1}$ (the previous formula has to be read as $\f(\e,\t)=c+\e w$ when $n=1$) are trivial solutions of \eqref{ESM}, and that they parametrize the so called ``vertical hyperplanes'', i.e. (left- translations of) maximal subgroups of $\heis^n$: it follows that these hypersurfaces are stationary points of the area functional, and a simple calibration argument implies that they are also minimizers since they have constant horizontal normal (cfr. Example \ref{exconstnorm}). These considerations suggest that the right counterpart of the classical Bernstein problem in the Heisenberg setting is\medskip

{\bfseries Bernstein problem for $X_1$-graphs in $\heis^n$:} are there global solutions $\f:\R^{2n}\to\R$ of the minimal surface equation \eqref{ESM} which cannot be written as in \eqref{iperpiani}? \medskip

As we will see, again the answer seems to depend on the dimension $n$ of the space; however, new and unexpected phenomena seem to arise, e.g. the fact that we have solutions to \eqref{ESMC} which are not area minimizing.

\subsection{The Bernstein problem in $\heis^1$}\label{bernstein1}
We have seen in Section \ref{solint} that for $n=1$ there exist solutions of \eqref{ESM} which cannot be written as in \eqref{iperpiani}; see for instance Examples \ref{solint1} and \ref{solint2}. 

We already pointed out that every solution of the classic minimal surface equation \eqref{ESMC} parametrizes (the boundary of) a globally minimizer; in $\heis^1$ instead a new phenomenon arise, in the sense that there are entire solutions of the intrinsic minimal surface equation \eqref{ESM} which parameterize a surface which is not a minimizer. However whenever the surface is $\H$- perimeter minimizing in $\H^1$ it has to be a vertical plane. More precisely, we have the following

\begin{theo}[Minimizers vs. stationary entire $X_1$-graphs]\label{bernsteinminim}
Let $\f:\R^2\to\R$ be a $\ci^2$ function and let $E,S\subset\H^1$ be respectively the $X_1$-graph and the $X_1$-subgraph induced by $\f$ as in (\ref{Xgraph}) and (\ref{Xsubgraph}). Let us suppose $E$ is a minimizer for the $\heis$-perimeter measure in $\H^1\equiv\R^3$, i.e. $E$ satisfies (\ref{minimizer}) with $\Omega=\R^3$.


Then $S$ is a vertical plane, i.e. $\f(\e,\t)=w\e+c$ for all $(\e,\t)\in\R^2$ for some constants $w,c\in\R$.
\end{theo}
\begin{proof}

STEP 1. First of all, we want to rewrite the second variation formula \eqref{varsec0.4} in the coordinates $c,t$ introduced in Section 3. Therefore let $G$ be defined by
\begin{eqnarray*}
G &:& \R^2_{x,t}\to\R^2_{\e,\t}\\
&& (x,t)\longmapsto(t,-4x)
\end{eqnarray*}
and set
$$A(x):=(\Wf\f\circ G)(x,0),\qquad B(x):=(\f\circ G)(x,0);$$
in particular, $\f\circ G$ is an entire solution of \eqref{DB}. As in Section \ref{solint} we set $x(c,t):=\frac{A(c)}{2}t^2+B(c)t+c$ and 
\begin{eqnarray*}
F &:& \R^2_{c,t}\to\R^2_{x,t}\\
&& (c,t)\longmapsto(x(c,t),t)
\end{eqnarray*}
 
Therefore, if we define
$$\begin{array}{l}
\Omega:=F(\R^2)\subset\R^2_{x,t},\\
F^\ast:=G\circ F,\\
\Omega^\ast:=F^\ast(\R^2)=G(\Omega)\subset\R^2_{\e,\t}
\end{array}$$
and $c:\Omega\to\R$ through the formula $F^{-1}(x,t)=(c(x,t),t)$, thanks to Theorem \ref{teonote1.1bis} we have
\begin{itemize}
\item for any $c\in\R$ we have $A'(c)=B'(c)=0$ or $B'(c)^2< 2A'(c)$;
\item $F^\ast$ is a $\ci^2$ diffeomorphism between $\R^2_{c,t}$ and $\Omega^\ast$. Moreover, $\Omega$ and $\Omega^\ast$ are open neighbourhood of the lines $\{t=0\}$ and $\{\e=0\}$ respectively.
\end{itemize}

It is not difficult to prove that for all $(\e,\t)\in\Omega^\ast$ one has
\begin{eqnarray}
&& \f(\e,\t)=A(c(-\t/4,\e))\e+B(c(-\t/4,\e)))=\frac{\partial x}{\partial t}(F^{\ast-1}(\e,\t));\label{8.4.1}\\
&& \Wf\f(\e,\t)=A(c(-\t/4,\e)).\label{8.4.2}
\end{eqnarray}
and taking into account that
\begin{eqnarray}
&& \frac{\partial c}{\partial x}(x,t)=\frac{1}{\frac{\partial x}{\partial c}(F^{-1}(c,t))}=\frac{1}{\frac{A'(c(x,t))}{2}t^2 + B'(c(x,t)) +1}\label{8.4.3}\\
&& \frac{\partial c}{\partial t}(x,t)=-\frac{\frac{\partial x}{\partial t}(F^{-1}(c,t))}{\frac{\partial x}{\partial c}(F^{-1}(c,t))}=-\frac{A(c(x,t))t+B(c(x,t))}{\frac{A'(c(x,t))}{2}t^2 + B'(c(x,t)) +1}\label{8.4.4}
\end{eqnarray}
for all $x,t\in\Omega$, we get for all $(\e,\t)\in\Omega^\ast$ that
\begin{eqnarray}
\Tt\Wf\f + 2(\Tt\f)^2 &=& -\frac{1}{4} \frac{A'(c(F^{\ast-1}(\e,\t)))}{\frac{\partial x}{\partial c}(F^{\ast-1}(\e,\t))}+ \nonumber\\
&& + 2\left[\frac{1}{4} \frac{A'(c(F^{\ast-1}(\e,\t))) + B'(c(F^{\ast-1}(\e,\t)))}{\frac{\partial x}{\partial c} (F^{\ast-1}(\e,\t))} \right]^2\nonumber\\
&=& \frac{1}{8} \frac{-2A'(c)\frac{\partial x}{\partial c}+(\frac{\partial^2 x}{\partial c\partial t})^2} {(\frac{\partial x}{\partial c})^2}(F^{\ast-1}(\e,\t)).\label{8.4.5}
\end{eqnarray}
Observe that for any $(c,t)\in\R^2$ we have
\begin{eqnarray}
\frac{1}{8}\frac{-2A'(c)\frac{\partial x}{\partial c}(c,t)+(\frac{\partial^2 x}{\partial c\partial t}(c,t))^2}{(\frac{\partial x}{\partial c}(c,t))^2} &=&  \frac{1}{8}\frac{-2A'(c)(\frac{A'(c)}{2}t^2+B'(c)t+1)+(A'(c)t+B'(c))^2}{(\frac{A'(c)}{2}t^2+B'(c)t+1)^2}\nonumber\\
&=& \frac{1}{8}\frac{-2A'(c)+B'(c)^2}{(\frac{A'(c)}{2}t^2+B'(c)t+1)^2}\ \leq 0 \label{8.4.6}
\end{eqnarray}

Notice that the correspondance 
$$\ci^1_c(\R^2_{c,t})\ni\z\longleftrightarrow\p:=\z\circ F^{\ast-1}\in\ci^1_c(\Omega^\ast)$$
is bijective and
\begin{eqnarray}
(\Wf\p)(F^\ast(c,t)) &=& \frac{\partial \p}{\partial\e}(F^\ast(c,t))-4\f(F^\ast(c,t))\frac{\partial \p}{\partial\t}(F^\ast(c,t))\nonumber\\
&=& \frac{\partial \p}{\partial\e}(F^\ast(c,t))-4\frac{\partial x}{\partial t}(c,t)\frac{\partial \p}{\partial\t}(F^\ast(c,t))\nonumber\\
&=& \frac{\partial\z}{\partial t}(c,t).\label{8.4.8}
\end{eqnarray}
Since
\begin{eqnarray*}
\text{det }JF^\ast(c,t) &=& \text{det }JG(F(c,t))\text{det }JF(c,t)\\
&=& 4(\textstyle\frac{A'(c)}{2}t^2+B'(c)t+1)\ >0
\end{eqnarray*}
a change of variable and equations \eqref{varsec0.4}, \eqref{8.4.5}, \eqref{8.4.6} and \eqref{8.4.8} give
\begin{eqnarray}
g''(0) &=& \int_{\Omega^\ast} \frac{(\Wf\p)^2 +8\p^2\:[\Tt\Wf\f+2(\Tt\f)^2]}{\left[1+|\Wf\f|^2\right]^{3/2}}\ d\e d\t\nonumber\\
&=& 4\int_{\R^2} \frac{(\frac{\partial\z}{\partial t})^2+\z^2\frac{-2A'(c)+B'(c)^2}{(\frac{A'(c)}{2}t^2+B'(c)t+1)^2}} {[1+A(c)^2]^{3/2}} \left[\frac{A'(c)}{2}t^2+B'(c)t+1\right]\ dc\,dt\nonumber\\
&=& 4\int_{\R^2}\left[\left(\frac{\partial\z}{\partial t}\right)^2 u + \z^2v \right]\ dc\,dt\label{8.4.9}
\end{eqnarray}
where $g$ is as in \eqref{g(s)} and we have set $\z:=\p\circ F^\ast$ and
$$u(c,t):=\frac{\frac{A'(c)}{2}t^2+B'(c)t+1}{[1+A(c)^2]^{3/2}}$$
$$v(c,t):=\frac{B'(c)^2-2A'(c)}{[1+A(c)^2]^{3/2}[\frac{A'(c)}{2}t^2+B'(c)t+1]}.$$

The fact that $\f$ parameterizes a minimizer implies that $g''(0)\geq 0$ for all $\p\in\ci^1_c(\Omega^\ast)$; since $F^\ast:\R^2\to\Omega^\ast$ is a $\ci^2$ diffeomorphism we deduce that
\begin{equation}\label{8.4.10}
\int_{\R^2}\left[\left(\frac{\partial\z}{\partial t}\right)^2 u + \z^2v \right]\ dc\,dt\geq 0 \qquad\forall\z\in\ci^1_c(\R^2).
\end{equation}

STEP 2. It is easy to see that our thesis on $\f$ is equivalent to $A$ and $B$ being constant, i.e. to $A'=B'\equiv 0$. Suppose by contradiction that there exist a $c_0\in\R$ such that this does not hold, then by Theorem \ref{teonote1.1bis} we have $b^2<2a$, where $b:=B'(c_0)$ and $a:=A'(c_0)>0$. We want to use the second variation formula \eqref{8.4.10} to obtain simpler conditions, namely inequalities on certain one-dimensional integrals involving $a$ and $b$ (see equation \eqref{8.4.15}).

Fix therefore a $\z\in\ci^1_c(\R^2)$ and set
$$\z_\ep(c,t):=\frac{1}{\sqrt{\ep}}\z\bigl( c_0+\textstyle\frac{c-c_0}{\ep},t\bigr);$$
then by \eqref{8.4.10} we get 
\begin{equation}\label{8.4.11}
0\leq \int_{\R^2}\left(\frac{\partial\z_\ep}{\partial t}\right)^2 u\ dc\,dt +\int_{\R^2} \z_\ep^2v \ dc\,dt =: I_\ep + II_\ep
\end{equation}

Observe that
\begin{eqnarray*}
I_\ep &=& \frac{1}{\ep} \int_{\R^2}\left(\frac{\partial\z}{\partial t}\bigl( c_0+\textstyle\frac{c-c_0}{\ep},t\bigr) \right)^2 u(c,t)\ dc\,dt\\
&=& \int_{\R^2} \left(\frac{\partial\z}{\partial t}(u,t) \right)^2 u(c_0+\ep(u-c_0),t)\ du\,dt
\end{eqnarray*}
and by Lebesgue convergence theorem one obtains
\begin{equation}\label{8.4.12}
\lim_{\ep\to 0}I_\ep = \int_{\R^2}\left(\frac{\partial\z}{\partial t}(c,t)\right)^2 u(c_0,t)\ dc\,dt.
\end{equation}
Analogously one gets
\begin{equation}\label{8.4.13}
\lim_{\ep\to 0}II_\ep = \int_{\R^2}z(c,t)^2 v(c_0,t)\ dc\,dt.
\end{equation}

Combining \eqref{8.4.11}, \eqref{8.4.12} and \eqref{8.4.13} we get
\begin{equation}\label{8.4.14}
\int_{\R^2} \left(\frac{\partial\z}{\partial t}(c,t)\right)^2 h(t)\ dc\,dt \geq (2a-b^2) \int_{\R^2}z(c,t)^2 \frac{1}{h(t)}\ dc\,dt
\end{equation}
for all $\z\in\ci^1_c(\R^2)$, where we have put
$$h(t):=\frac{a}{2}t^2+bt+1.$$
By standard arguments (taking for example $\z(c,t)$ of the form $\z_1(c)\z_2(t)$) we can infer the one-dimensional inequalities
\begin{equation}\label{8.4.15}
\int_\R \z'^2 h\ dt\geq (2a-b^2)\int_\R \z^2\frac{1}{h}\ dt\qquad\text{for all }\z\in\ci^1_c(\R).
\end{equation}

STEP 3. We will follow here the technique used in \cite{DGN2006} to provide a counterexample to \eqref{8.4.15}, which will give a contradiction.

For $\ep>0$ fix $\chi_\ep\in\ci^1_c(\R)$ such that
$$\begin{array}{l}
0\leq\chi_\ep\leq 1\\
\chi_\ep\equiv 1\text{ on }\bigl(-\frac{1}{\ep},\frac{1}{\ep}\bigr),\quad \text{spt } \chi_\ep\Subset \bigl(-\frac{2}{\ep},\frac{2}{\ep}\bigr)\\
|\chi_\ep'|\leq C\ep,\quad C>0\text{ independent of }\ep
\end{array}$$
and set 
$$\z_\ep(t):=\frac{\chi_\ep(t)}{\sqrt{h(t)}}.$$

Equation \eqref{8.4.15} becomes then
\begin{equation}\label{8.4.17}
\int_\R \z_\ep'^2 h\ dt\geq (2a-b^2)\int_\R \z_\ep^2\frac{1}{h}\ dt
\end{equation}
and observe that
\begin{equation}\label{8.4.18}
\lim_{\ep\to 0}\int_\R \z_\ep^2\frac{1}{h}\ dt = \int_\R\frac{dt}{\bigl(\frac{a}{2}t^2+bt+1\bigr)^2}.
\end{equation}

As for the left hand side of \eqref{8.4.17}, we have
\begin{equation}\label{8.4.19}
\int_\R\z_\ep'^2h\,dt=\int_\R\left( \frac{\chi_\ep'}{\sqrt{h}}-\frac{\chi_\ep h'}{2h^{3/2}}\right)^2h\,dt=\int_\R \chi_\ep'^2\,dt-\int_\R\frac{\chi_\ep\chi_\ep'h'}{h}\,dt +\frac{1}{4}\int_\R\chi_\ep^2\frac{h'^2}{h^2}\,dt;
\end{equation}
an integration by parts gives 
$$\int_\R\chi_\ep\chi_\ep'\frac{h'}{h}\,dt= - \frac{1}{2}\int_\R\chi_\ep^2\frac{h''}{h} \,dt+\frac{1}{2}\int_\R\chi_\ep^2\frac{h'^2}{h^2} \,dt$$
whence \eqref{8.4.19} rewrites as
$$\int_\R\z_\ep'^2h\,dt=\int_\R \chi_\ep'^2\,dt +\frac{1}{2}\int_\R\chi_\ep^2\frac{h''}{h} \,dt-\frac{1}{4}\int_\R\chi_\ep^2\frac{h'^2}{h^2}\,dt$$
Finally, by Lebesgue convergence theorem we infer
\begin{equation}\label{8.4.20}
\lim_{\ep\to 0} \int_\R\z_\ep'^2h \,dt= \frac{1}{2}\int_\R\frac{h''}{h} \,dt-\frac{1}{4}\int_\R\frac{h'^2}{h^2}\,dt = \frac{1}{4}\int_\R\frac{h''}{h}\,dt
\end{equation}
where, in the last equality, we integrated by parts again.

From \eqref{8.4.17}, \eqref{8.4.18} and \eqref{8.4.20} we obtain therefore
\begin{equation}\label{8.4.21}
\frac{1}{4}\int_\R\frac{a\,dt}{\frac{a}{2}t^2+bt+1} \geq (2a-b^2)\int_\R \frac{dt}{\bigl( \frac{a}{2}t^2+bt+1\bigr)^2}\,.
\end{equation} 
Since for $\a>0$ we have
$$\int_\R\frac{dt}{1+\a t^2}=\frac{\pi}{\sqrt{\a}}\quad\text{and}\quad \int_\R\frac{dt}{(1+\a t^2)^2}=\frac{\pi}{2\sqrt{\a}}$$
and observing that
$$\int_\R \frac{dt}{\bigl( \frac{a}{2}t^2+bt+1\bigr)^m}=\left(\frac{2a}{2a-b^2}\right)^m\int_\R\frac{dt}{(1+\a t^2)^m}\qquad m=1,2$$
with $\a:=\frac{a^2}{2a-b^2}$, by \eqref{8.4.21} we obtain
$$\frac{a}{4}\frac{2a}{2a-b^2}\pi\frac{\sqrt{2a-b^2}}{a}\geq(2a-b^2)\frac{4a^2}{(2a-b^2)^2}\frac{\pi}{2}\frac{\sqrt{2a-b^2}}{a}$$
which reduces to $1/2\geq 2$ recalling that $a>0$ and then a contradiction.

STEP 4. We have proved that $A$ and $B$ are constant functions, and this in turn implies that $\Omega^\ast=\R^2$ and $\phi(\e,\t)=A\e+B$. This completes the proof of the Theorem.
\end{proof}

%

\subsection{The Bernstein problem in $\heis^n$ for $n\geq 2$}
Let us exploit equation \eqref{ESM} and write it as
\begin{multline}\label{ESMespl}
\sum_{j=2}^n \Xt_j \left( \frac{\Xt_j\f}{\sqrt{1+|\Wf\f|^2}}\right) + \Wf_{n+1}\left(\frac{\Wf_{n+1}\f}{\sqrt{1+|\Wf\f|^2}}\right) +\\ 
+ \sum_{j=2}^n \Yt_j \left( \frac{\Yt_j\f}{\sqrt{1+|\Wf\f|^2}}\right)=0\qquad
\end{multline}
where $\f:\R^{2n}=\R_\e\times\R_v^{2n-2}\times\R_\t\to\R$ is of class $\ci^2$. Notice that, if one looks for solutions $\f$ which do not depend on the $\t$ variable, i.e. such that $\f(\e,v,\t)=\p(\e,v)$ for some $\p:\R^{2n-1}\to\R$, equation \eqref{ESMespl} rewrites as the classic minimal surface equation
\begin{equation}
\text{div}\left(\frac{\nabla\p}{\sqrt{1+|\nabla\p|^2}}\right)=0.
\end{equation}

This observation allows us to easily construct a counterexample to the Bernstein problem for $X_1$-graphs in $\heis^n$ when $n\geq 5$; in fact in this case we have $2n-1\geq 9$ and Theorem \ref{bernsteincl} provides a function $\p:\R^{2n-1}\to\R$ which solves \eqref{ESMC} and is not affine, i.e. the related $\f(\e,v,\t)=\p(\e,v)$ solves \eqref{ESMespl} and cannot be written as in \eqref{iperpiani}. 

We also notice that $X_1$-graphs of functions $\f(\e,v,\t)=\p(\e,v)$ (where $\p$ solves \eqref{ESMC}) are actually minimizers of the $\heis$-perimeter; in fact it is easy to check that the smooth section $\nu:\heis^n\to H\heis^n$ defined by
\begin{eqnarray*}
\nu(x,y,t) &=& \left( -\frac{1}{\sqrt{1+|\Wf \f|^2}},\frac{\Wf\f}{\sqrt{1+|\Wf \f|^2}} \right)(\e,v,0)\\
&=& \left(-\frac{1}{\sqrt{1+|\nabla\p|^2}},\frac{\nabla\p}{\sqrt{1+|\nabla\p|^2}} \right)(\e,v),
\end{eqnarray*}
where we put $\e:=y_1$ and $v:=(x_2,\dots,x_n,y_2,\dots,y_n)$, is a calibration for the graph of $\f$, i.e.
\begin{itemize}
\item $\div_X\ \nu=0$;
\item $|\nu(P)|=1$ for all $P\in\heis^n$;
\item $\nu$ coincides with the horizontal inward normal to the $X_1$-graph of $\f$ (see Theorem \ref{teoASCV}).
\end{itemize}
Observe that in this argument (which is basicly the same used to prove the minimality of any entire graph solution of \eqref{ESMC} in the classical case) it was essential the non-dependance of $\f$ on the vertical variable $\t$: as we have seen in Section \ref{bernstein1}, in general it is not true that an entire solution of \eqref{ESM} parameterizes a minimizer.

The Bernstein problem for intrinsic graphs in $\heis^n$, as far as we know, remains unsolved for $n=2,3,4$; observe that any possible negative answer must effectively depend on the variable $\t$, or the previous argument leading to the classic Bernstein equation could apply, contradicting Theorem \ref{bernsteincl}. 

\begin{obs}\label{prodottigruppo} 
{\upshape We explicitely notice that, as $\heis^n$ cannot be written as a product of $\heis^m$ for any $m<n$, in general there is no easy way to construct a solution to \eqref{ESM} in $\heis^n$ ``modifying'' somehow a given solution to \eqref{ESM} in $\heis^m$, e.g. managing the nontrivial ones found in Subsection \ref{bernstein1} for $m=1$.}
\end{obs}.

\end{document}